\documentclass[11pt,reqno]{amsart}
\usepackage{amsthm,amssymb,amsmath}
\usepackage{textcomp}
\usepackage{xcolor}

\usepackage[T1]{fontenc}
\usepackage{lmodern}

\newtheorem{theorem}{Theorem}
\newtheorem*{theorem*}{Theorem}
\newtheorem{lemma}{Lemma}
\newtheorem{corollary}{Corollary}
\newtheorem{proposition}{Proposition}

\theoremstyle{definition}

\newtheorem*{definition*}{\bf Definition}

\newtheorem{remark}{\sc Remark}
\newtheorem*{remark*}{\sc Remark}
\newtheorem*{remarks}{\sc Remarks}

\newtheorem*{example*}{\bf Example}

\newcommand{\loc}{{\rm loc}}

\newcommand{\dist}{\mbox{dist}}

\begin{document}

\title{Parabolic equations and SDEs with time-inhomogeneous Morrey drift}

\begin{abstract}
We prove the unique weak solvability of stochastic differential equations with time-inhomogeneous drift in essentially the largest (scaling-invariant) Morrey class, i.e.\,with integrability parameter $q>1$ close to $1$. The constructed weak solutions constitute a Feller evolution family. The proofs are based on a detailed Sobolev regularity theory of the corresponding parabolic equation.
\end{abstract}

\author{D.\,Kinzebulatov}

\email{damir.kinzebulatov@mat.ulaval.ca}

\address{Universit\'{e} Laval, D\'{e}partement de math\'{e}matiques et de statistique, Qu\'{e}bec, QC, Canada}

\keywords{Stochastic differential equations, weak solutions, singular drifts, Morrey class, Feller property}

\subjclass[2010]{60H10, 47D07 (primary), 35J75 (secondary)}

\thanks{The research of the author is supported by the NSERC (grant RGPIN-2017-05567).}

\maketitle

\section{Introduction}

\medskip

\textbf{1.~}We consider the problem of weak well-posedness of stochastic differential equation 
\begin{equation}
\label{sde1}
X_t=x-\int_0^t b(r,X_r)dr + \sqrt{2}B_t, \quad x \in \mathbb R^d,
\end{equation} 
where $B_t$ is a Brownian motion in $\mathbb R^d$,
under minimal assumptions on the time-inhomogeneous vector field $b:\mathbb R \times \mathbb R^d \rightarrow \mathbb R^d$ (drift), $d \geq 3$. 
This equation or, more generally, stochastic equations additionally having variable, possibly discontinuous  diffusion coefficients, arise e.g.\,in the problems of stochastic optimization and serve as a basis for many physical models. This requires, generally speaking, 
dealing with irregular, locally unbounded drifts.
An illustrative example is equation \eqref{sde1} with velocity field $b$ obtained by solving 3D Navier-Stokes equations, which models the motion of a small particle in a turbulent flow \cite{MK}. 
One is thus led to the problem of establishing weak and strong well-posedness of \eqref{sde1} under minimal assumptions on $b$. The latter can also be stated as the problem of finding the most general integral characteristics of $b$ that determine whether \eqref{sde1} is weakly/strongly well-posed. 

Let us give a brief outline of the literature on stochastic differential equations (SDEs) with singular drift. 
We will try to keep the chronological order, but will be somewhat loose with the terminology by including in ``well-posedness''  the uniqueness results of different strength.

The ``sub-critical'' Ladyzhenskaya-Prodi-Serrin class
\begin{equation}
\label{L}
|b| \in L^l(\mathbb R,L^p(\mathbb R^d)), \quad p \geq d, l \geq 2, \quad \frac{d}{p}+\frac{2}{l} < 1
\end{equation}
was attained by Portenko \cite{P} (weak solutions) and Krylov-R\"{o}ckner \cite{KR} (strong solutions).  See also Zhang \cite{Z,Z0,Z2}. 
Between \cite{P} and \cite{KR}, Bass-Chen \cite{BC} proved existence and uniqueness in law of weak solutions of \eqref{sde1} for $b=b(x)$ in the Kato class of vector fields.
The Kato class contains $\{|b| \in L^p(\mathbb R^d), p>d\}$ as well as some vector fields with entries not even in $L^{1+\varepsilon}_{\loc}(\mathbb R^d)$, $\varepsilon>0$, however, it does not contain  $\{|b| \in L^d(\mathbb R^d)\}$. Speaking of time-homogeneous drifts, of course, the fact that $p=d$ is the optimal exponent on the scale of Lebesgue spaces can be seen from rescaling the parabolic equation. 

In \cite{BFGM}, Beck-Flandoli-Gubinelli-Maurelli  developed an approach to proving strong well-posedness of \eqref{sde1} with $b$ in the critical Ladyzhenskaya-Prodi-Serrin class
\begin{equation}
\label{L2}
|b| \in L^l(\mathbb R,L^p(\mathbb R^d)), \quad p \geq d, l \geq 2, \quad \frac{d}{p}+\frac{2}{l} \leq 1.
\end{equation}
 for a.e.\,starting point $x \in \mathbb R^d$
via stochastic transport and stochastic continuity equations. They also considered the following example. Let
\begin{equation}
\label{b_hardy}
b(x)=\sqrt{\delta} \frac{d-2}{2}\mathbf{1}_{|x|<1}|x|^{-2}x,
\end{equation}
so $|b|$ just misses to be $L^d(\mathbb R^d)$.
If $\delta>4(\frac{d}{d-2})^2$, i.e.\,the attraction to the origin by the drift is large enough, then SDE \eqref{sde1} with the starting point $x=0$ does not have a weak solution.
In \cite{KiS1}, Sem\"{e}nov and the author showed that the Feller generator $\Delta - b \cdot \nabla$ with ``weakly form-bounded'' $b=b(x)$, see \eqref{wfbd} below, determines for every starting point $x \in \mathbb R^d$ a weak solution to \eqref{sde1} that is, moreover, unique among weak solutions that can be constructed via approximation. (To the best of the author's knowledge, this was the first result on weak well-posedness of \eqref{sde1} that included both $|b| \in L^d(\mathbb R^d)$ and the model vector field \eqref{b_hardy} with $\delta$ small; it also included the elliptic Morrey class with $q>1$ and the Kato class considered by Bass-Chen.) 
Returning to time-inhomogeneous drifts, we note that almost at the same time Wei-Lv-Wu \cite{WLW}, and later Nam \cite{N}, obtained results on weak well-posedness of \eqref{sde1} for every $x \in \mathbb R^d$ for time-inhomogeneous vector fields $b$ that can be more singular than the ones in \eqref{L}. Nevertheless, their results excluded $b=b(x)$ with $|b| \in L^d(\mathbb R^d)$. 
In \cite{Kr1_}, Krylov proposed a new approach to proving strong well-posedness of \eqref{sde1} e.g.\,with $|b| \in L^d(\mathbb R^d)$ based on his and Veretennikov's old criterion for a weak solution to be a strong solution  \cite{VK}.
In \cite{XXZZ}, Xia-Xie-Zhang-Zhao established, among other results, weak well-posedness of \eqref{sde1} for every initial point and $b \in C_b(\mathbb R,L^d(\mathbb R^d))$. R\"{o}ckner-Zhao \cite{RZ} furthermore established weak well-posedness of \eqref{sde1}, with any $x \in \mathbb R^d$, for drifts in $L^\infty(\mathbb R,L^{d,w}(\mathbb R^d))$, plus the drifts in the critical LPS class. Here $L^{d,w}(\mathbb R^d)$ denotes the weak $L^d$ class that contains vector fields in $L^d(\mathbb R^d)$, as well as more singular vector fields, such as \eqref{b_hardy}. In \cite{RZ2}, they obtained strong well-posedness of \eqref{sde1}, for any starting point and $b$ in the critical LPS class. In \cite{KiM}, the author and Madou established weak well-posedness of \eqref{sde1}, for every starting point and form-bounded drifts, see example 5) below. This class contains $L^\infty(\mathbb R,L^{d,w}(\mathbb R^d))$  as well as some drifts that are not even in $L^\infty(\mathbb R,L^{2+\varepsilon}(\mathbb R^d))$ for a given $\varepsilon>0$.
By the way, in \cite{KiSS}  Sem\"{e}nov, Song and the author showed that the approach of \cite{BFGM} to the strong well-posedness of \eqref{sde1} via the stochastic transport/continuity equations also works for form-bounded vector fields $b=b(x)$, although, again, one obtains strong well-posedness of \eqref{sde1} only for a.e.\,starting point.

The above outline does not discuss many interesting results on distribution-valued drifts and  drifts satisfying additional assumptions on their structure, such as ${\rm div\,}b \leq 0$. The latter include results of partial well-posedness of \eqref{sde1} with $|b| \in L^p(\mathbb R^d)$ in the supercritical regime (in the sense of scaling) $\frac{d}{2}<p<d$. In this regard, see Zhang-Zhao \cite{ZZ} and Zhao \cite{Zhao}.

In the present paper we consider drifts in the Morrey class $E_q$ with integrability parameter $q>1$ that can be chosen arbitrarily close to $1$. 
Denote
$$
C_r(t,x):=\{(s,y) \in \mathbb R^{d+1} \mid t \leq s \leq t+r^2, |x-y| \leq r\}
$$
and, given a vector field $b:\mathbb R^{d+1} \rightarrow \mathbb R^d$ with components in $L^q_{\loc}(\mathbb R^{d+1})$,
 $q \in [1,d+2]$, set
\begin{align*}
\|b\|_{E_q}&:=\sup_{r>0, z \in \mathbb R^{d+1}} r\biggl(\frac{1}{|C_r|}\int_{C_r(z)}|b(t,x)|^q dtdx \biggr)^{\frac{1}{q}}\\
&=\sup_{r>0, z \in \mathbb R^{d+1}} r\biggl(\frac{1}{|C_r|}\int_{C_r(z)}|b(-t,x)|^q dtdx \biggr)^{\frac{1}{q}}.
\end{align*}

\begin{definition*}
We say that a vector field $b$ belongs to the parabolic Campanato-Morrey (or, for brevity, Morrey) class $E_{q}$ if
$
\|b\|_{E_q}<\infty.
$
\end{definition*}

One has $$\|b\|_{E_{q}} \leq \|b\|_{E_{q_1}} \text{ if } q<q_1.$$

If above $b=b(x)$, then one obtains the usual elliptic Morrey class $M_q$, that is, $|b| \in L^q_{\loc}(\mathbb R^d)$ and
$$
\|b\|_{M_q}:=\sup_{r>0, y \in \mathbb R^{d}} r\biggl(\frac{1}{|B_r|}\int_{B_r(y)}|b(x)|^q dx \biggr)^{\frac{1}{q}}<\infty,
$$
where $B_r(y)$ is the closed ball of radius $r$ centered at $y$.

Our result, stated briefly, is as follows (see Theorems \ref{thm1}-\ref{thm3} for details). 
We will be using some notations defined in the end of this section.

\begin{theorem*}
Let $d \geq 3$, let $b:\mathbb R^{d+1} \rightarrow \mathbb R^d$ be a vector field in the Morrey class $E_q$ with $q>1$ close to $1$. Let $p \in ]1,\infty[$. There exists a constant $c_{d,p,q}$ such that if $\|b\|_{E_q}<c_{d,p,q}$, then the following are true:

{\rm(\textit{i})} There exists a unique weak solution  to 
\begin{align*}
(\lambda - \partial_t - \Delta + b(t,x) \cdot \nabla)u&=0, \quad t<r, \\
 u(r,\cdot)&=g(\cdot) \in L^p(\mathbb R^d) \cap L^2(\mathbb R^d).
\end{align*}
The difference
$$
u(t,\cdot)-\frac{e^{-\lambda(r-t)}}{(4 \pi (r-t))^{\frac{d}{2}}}\int_{\mathbb R^d} e^{-\frac{|\cdot-y|^2}{4(r-t)}}g(y)dy \quad (t<r), \quad \text{extended by $0$ to $t>r$},
$$
 is in the parabolic Bessel potential space $ \mathbb W^{1+\frac{1}{p},p}(\mathbb R^{d+1})$.

{\rm(\textit{ii})} For $p>d+1$, operators $\{P^{t,r}\}_{t<r}$ defined by $$P^{t,r}g:=u(t), \quad g \in C_\infty(\mathbb R^d) \cap L^2(\mathbb R^d)$$ are extended to a backward Feller evolution family on $C_\infty(\mathbb R^d)$ that determines, for every $x \in \mathbb R^d$, a weak solution to SDE \eqref{sde1}
that is, moreover, unique in some large classes of weak solutions. 
\end{theorem*}

Here are some examples of  vector fields in $E_q$, $q>1$:

\smallskip

1) The critical Ladyzhenskaya-Prodi-Serrin class
$$
|b| \in L^l(\mathbb R,L^p(\mathbb R^d)), \quad p \geq d, l \geq 2, \quad \frac{d}{p}+\frac{2}{l} \leq 1.
$$
To prove the inclusion of this class into $E_q$ it suffices to consider, by an elementary interpolation argument, only the cases $l=2$, $p=\infty$ and $l=\infty$, $p=d$.
In the former case the inclusion is trivial, in the latter case the inclusion follows using H\"{o}lder's inequality.

This example is strengthened in the next two examples.

\smallskip

2) The vector fields $b$ with $|b| \in L^{2,w}(\mathbb R,L^\infty(\mathbb R^d))$ are in $E_q$ for $1<q<2$.  Indeed, by a well known characterization of weak Lebesgue spaces, we have
\begin{align*}
r \biggl(\frac{1}{|C_r|}\int_{C_r}|b|^q dz \biggr)^{\frac{1}{q}} & \leq C r\biggl(\frac{1}{r^{2}}\int_t^{t+r^2}|\tilde{b}|^q ds \biggr)^{\frac{1}{q}}  \qquad \tilde{b}(t):=\|b(t,\cdot)\|_{L^\infty(\mathbb R^d)}\\
&\leq C \|\tilde{b}\|_{L^{2,w}(\mathbb R)}.
\end{align*}
Hence, for example, a vector field $b$ that satisfies
$$
\|b(t,\cdot)\|_{L^\infty(\mathbb R^d)} \leq \frac{C}{\sqrt{t}}, \quad t>0
$$
 (and defined to be zero for $t \leq 0$)
is in $E_q$ with $1<q<2$.

\smallskip

3) Moreover, by the well known inclusion of the weak Lebesgue space $L^{d,w}(\mathbb R^d)$ in the elliptic Morrey class, $$|b| \in L^\infty(\mathbb R,L^{d,w}(\mathbb R^d)) \quad \Rightarrow \quad b \in E_q\;\text{ with }1<q \leq d.$$

\smallskip

4) For every $\varepsilon>0$, one can find a vector field $b \in E_q$ such that $|b|$ is not in $ L^{q+\varepsilon}_{\loc}(\mathbb R^{d+1})$. (In particular, selecting $q>1$ sufficiently close to $1$, we obtain vector fields $b$ satisfying the assumptions of the theorem, and that are not in $L^{1+\varepsilon}_{\loc}(\mathbb R^{d+1})$ for a given $\varepsilon>0$.)

\smallskip

5) A vector field $b$ is said to be form-bounded if $|b| \in L^2_{\loc}(\mathbb R^{d+1})$ and for a.e.\,$t \in \mathbb R$
$$
\|b(t,\cdot)\varphi\|_{L^2(\mathbb R^d)}^2 \leq \delta\|\nabla \varphi\|_{L^2(\mathbb R^d)}^2 + g(t) \|\varphi\|_{L^2(\mathbb R^d)}^2
$$ 
for all $\varphi=\varphi(\cdot) \in C_c^\infty(\mathbb R^d)$,
for some constant $\delta>0$ and a function $g \in L^1_{\loc}(\mathbb R)$ (written as $b \in \mathbf{F}_\delta$). The constant $\delta$ is called a form-bound of $b$. 

This class itself contains 1) and 3). In particular, by Hardy's inequality,
the vector field \eqref{b_hardy} is in $\mathbf{F}_\delta$ with $g=0$ (but not in any $\mathbf{F}_{\delta'}$ with $\delta'<\delta$). 

Note that if a form-bounded vector field $b$ depends only on time, then $b \in L^2_{\loc}(\mathbb R)$. One can compare this with example 2).
\medskip

\textbf{2.~}Let us say a few more words about the form-bounded vector fields. In the time-homogeneous case $b=b(x)$ the form-boundedness of $b$
ensures, by the Lax-Milgram Theorem, that, whenever $\delta<1$, the formal operator $\Delta - b \cdot \nabla$ has a realization as the generator of a quasi-contraction semigroup in $L^2(\mathbb R^d)$ which, moreover, is bounded as an operator $W^{1,2}(\mathbb R^d) \rightarrow W^{-1,2}(\mathbb R^d)$. The form-boundedness of $b$ \textit{with} $\delta<1$ is essentially the broadest
condition on $|b|$ that provides the minimal ``classical'' theory of operator $-\Delta + b \cdot \nabla$ in $L^2$.

By considering dilates and translates of a  bump function, one can show that the class of form-bounded vector fields $\mathbf{F}_\delta$ is contained in $E_2$. (Let us note, on the other hand, that the time-homogeneous vector fields in $E_q$, $q>2$ are form-bounded \cite{F}, see also \cite{CF}. This remark extends right away e.g.\,to time-inhomogeneous vector fields $b=b(t,\cdot)$ having uniformly bounded in $t$ elliptic Morrey norm $\|\cdot\|_{M_q}$, $q>2$; they are, obviously, contained in $E_q$.) Note also that if $b$ is form-bounded with $g=0$, then  $\|b\|_{E_2}=c\sqrt{\delta}$ for appropriate constant $c=c_{d}$.

The class of form-bounded vector fields is well known in the literature on regularity theory of parabolic equations \cite{KS} (although, perhaps, less so than the Morrey class).
As we mentioned earlier, in the context of SDEs, condition $b \in \mathbf{F}_\delta$ with $\delta<d^{-2}$ was used to develop a Sobolev regularity theory of the corresponding parabolic equation and construct a Feller evolution family \cite{Ki}
that determines, for every $x \in \mathbb R^d$, a weak solution to \eqref{sde1} which is, moreover, unique in a broad class of weak solutions \cite{KiM}. 

\smallskip

\textbf{3.~}In the present paper 
we deal with time-inhomogeneous vector fields that can be more singular than the ones covered by the form-boundedness condition. These new vector fields are situated in $E_{q} - E_{s}$, $1<q\leq 2$, $s > 2$. 
Let us note that the distinction between $E_q$, $1<q < 2$ and smaller class $E_{s}$, $s > 2$ is not only quantitative.
Namely, consider a vector field $b=b(x)$ in $E_s$, $s>2$. Since it is form-bounded, the term $b \cdot \nabla$ in the parabolic equation can be handled using quadratic inequality: 
\begin{equation}
\label{q_ineq}
\int_{\mathbb R^{d}} (b \cdot \nabla u) u dx  \leq \gamma\int_{\mathbb R^{d}} |b|^2 u^2 dx  + \frac{1}{4\gamma}\int_{\mathbb R^{d}} |\nabla u|^2 dx, \quad \gamma>0,
\end{equation}
where one  estimates  $\int_{\mathbb R^{d}} |b|^2 u^2 dx$ from above, using form-boundedness of $b$, by $\int_{\mathbb R^{d}} |\nabla u|^2 dx$. This elementary argument plays a key role e.g.\,in the proof of gradient estimates in \cite{BFGM,Kr1_,Ki}, 
needed to prove well-posedness of the corresponding stochastic equations. However, quadratic inequality \eqref{q_ineq} is unavailable under our assumption $b \in E_q$, $1<q< 2$.

\smallskip

\textbf{4.~}The proofs in \cite{KiM} use a Feller evolution family constructed in \cite{Ki} using a parabolic variant of the iteration procedure of \cite{KS}. In this paper we pursue a  simpler operator-theoretic approach, replacing the iterations with the Duhamel series (written in ``resolvent form''). Namely, we construct solution to inhomogeneous equation 
$$
(\lambda + \partial_t - \Delta + b(t,x) \cdot \nabla)u=f \in L^p(\mathbb R^{d+1})
$$
(which is essential to the rest of the paper) as
\begin{equation}
\label{u_repr0}
u:=(\lambda+\partial_t-\Delta)^{-1}f - (\lambda+\partial_t-\Delta)^{-\frac{1}{2}-\frac{1}{2p}}Q_p (1+T_p)^{-1}R_p (\lambda+\partial_t-\Delta)^{-\frac{1}{2p'}} f,
\end{equation}
where $p':=\frac{p}{p-1}$ and, if $b \in E_q$, $q>1$, the operators
$$R_p=b^{\frac{1}{p}}\cdot \nabla (\lambda+\partial_t-\Delta)^{-\frac{1}{2}-\frac{1}{2p}}, \quad Q_p =(\lambda+\partial_t-\Delta)^{-\frac{1}{2p'}}|b|^{\frac{1}{p'}}$$ are bounded  on $L^{p}(\mathbb R^{d+1})$ and, provided $\|b\|_{E_q}$ is sufficiently small and $\lambda$ is large, the operator $T_p=R_p Q_p$ has norm $\|T_p\|_{p \rightarrow p} <1$, so the Duhamel series converges. The regularizing factor $(\lambda+\partial_t-\Delta)^{-1/2-1/2p}$ in \eqref{u_repr0} now yields the sought regularity of solution $u$ provided that $p$ is large ($>d+1$). A similar argument was used in \cite{Ki2} in the elliptic setting, where an even larger than $\{b=b(x)\} \cap E_q$, $q>1$ class of time-homogeneous vector fields was treated (see \eqref{wfbd} below); the weak well-posedness of SDEs with such drifts was addressed in \cite{KiS1}. The proof of the gradient estimates in \cite{Ki2}, however, depends on the symmetry of the resolvents, and the transition from elliptic estimates to the results for the parabolic equations required development of some new approaches to constructing semigroups, using old ideas of Hille (pseudoresolvents) and Trotter, see also \cite{KiS0}.  In the present paper the proofs are much shorter since we work directly with the parabolic operator. See also discussion after Theorem \ref{thm3}.

The necessity of the assumption ``$\|b\|_{E_q}$ cannot be too large'' follows from the aforementioned counterexample  to solvability of \eqref{sde1} with drift \eqref{b_hardy} when $x=0$ and $\delta>4(\frac{d}{d-2})^2$ (note that there $\|b\|_{E_q}=c\sqrt{\delta}$).
It was recently proved in \cite{KiS5} that, given an arbitrary $b \in \mathbf{F}_\delta$ with $\delta<4$,  SDE \eqref{sde1} has a weak solution for every starting point $x \in \mathbb R^d$, so the above example and this result become essentially sharp in high dimensions. (The fact that $\delta=4$ is critical can be seen from multiplying the parabolic equation \eqref{eq1} corresponding to \eqref{sde1} by $u|u|^{p-2}$, integrating by parts and using quadratic inequality and form-boundedness as above. The admissible $p$ that give e.g.\,an energy inequality turn out to be $p>\frac{2}{2-\sqrt{\delta}}$.)

\subsection*{Notations}

Set for $0<\alpha \leq 2$
\begin{equation}
\label{par_op_def0}
(\lambda-\partial_t-\Delta)^{-\frac{\alpha}{2}}h(t,x):=\int_{t}^\infty\int_{\mathbb R^d}e^{-\lambda(s-t)}\frac{1}{(4\pi (s-t))^{\frac{d}{2}}}\frac{1}{(s-t)^{\frac{2-\alpha}{2}}}e^{-\frac{|x-y|^2}{4(s-t)}}h(s,y)dsdy,
\end{equation}
\begin{equation}
\label{par_op_def}
(\lambda+\partial_t-\Delta)^{-\frac{\alpha}{2}}h(t,x):=\int_{-\infty}^t\int_{\mathbb R^d} e^{-\lambda(t-s)}\frac{1}{(4\pi (t-s))^{\frac{d}{2}}}\frac{1}{(t-s)^{\frac{2-\alpha}{2}}}e^{-\frac{|x-y|^2}{4(t-s)}}h(s,y)dsdy,
\end{equation}
where $\lambda \geq 0$. By a standard result, if $\lambda>0$, then these operators are bounded on $L^p(\mathbb R^{d+1})$, $1 \leq p \leq \infty$, with operator norm $\lambda^{-\frac{\alpha}{2}}$. If $\lambda>0$, then 
$(\lambda \pm \partial_t-\Delta)^{-1}$ is the resolvent of a Markov generator on $L^p(\mathbb R^{d+1})$, $1 \leq p<\infty$, which we will denote by $\lambda \pm \partial_t-\Delta$, respectively. (The abuse of notation resulting from not indicating $p$ should not cause any confusion.) In particular, one has well defined fractional powers $(\lambda \pm \partial_t-\Delta)^\frac{\alpha}{2}$. We refer to articles \cite{B,G}, among others, regarding the properties of these operators.

Denote by $\langle\, ,\rangle$ the integration in $d+1$ variables, i.e.
$$
\langle h \rangle :=\int_{\mathbb R^{d+1}} h dz, \quad \langle h,g\rangle:=\langle hg\rangle
$$
(all functions considered in this paper are real-valued).

We denote by $\mathcal B(X,Y)$ the space of bounded linear operators between Banach spaces $X \rightarrow Y$, endowed with the operator norm $\|\cdot\|_{X \rightarrow Y}$. Set  $\mathcal B(X):=\mathcal B(X,X)$.

Put $$\|h\|_p:=\langle |h|^p\rangle^{\frac{1}{p}}.$$ 
Denote by $\|\cdot\|_{p \rightarrow q}$ the $(L^p(\mathbb R^{d+1}),\|\cdot\|_p) \rightarrow (L^q(\mathbb R^{d+1}),\|\cdot\|_q)$ operator norm.

We write $T=s\mbox{-} X \mbox{-}\lim_n T_n$ for $T$, $T_n \in \mathcal B(X)$ if $Tf=\lim_n T_nf$ in $X$ for every $f \in X$.

Let $\lambda>0$ be fixed.
Set 
$$
\mathbb W^{\alpha,p}(\mathbb R^{d+1}):=(\lambda+\partial_t-\Delta)^{-\frac{\alpha}{2}}L^p(\mathbb R^{d+1})
$$
endowed with the norm $\|h\|_{\mathbb W^{\alpha,p}}:=\|(\lambda+\partial_t-\Delta)^{-\frac{\alpha}{2}}h\|_p$.

Denote by $C_\infty(\mathbb R^{d})$ the space of continuous functions on $\mathbb R^d$ vanishing at infinity, endowed with the $\sup$-norm. Define in the same way $C_\infty(\mathbb R^{d+1})$.

We fix $T>0$ and put 
$$D_T:=\{(s,t) \in \mathbb R^2 \mid 0 \leq s \leq t \leq T\}.$$

Recall that a family of positivity preserving $L^\infty$ contractions $\{U^{t,s}\}_{(s,t) \in D_T} \subset \mathcal B(C_\infty(\mathbb R^d))$ is called a Feller evolution family (on $C_\infty(\mathbb R^d)$) if $U^{t,r}U^{r,s}=U^{t,s}$ for all $r \in [s,t]$, $U^{s,s}={\rm Id}$ and 
$$
U^{r,s}=s\mbox{-}C_\infty(\mathbb R^d)\mbox{-}\lim_{t \downarrow r}U^{t,s}
$$
for all $s \leq r<T$.

\medskip

Define the following regularization of $b:\mathbb R^{d+1} \rightarrow \mathbb R^d$, $b \in E_q$, $q>1$:
\begin{equation}
\label{b_n}
b_n:=\mathbf{1}_nb, \quad \text{ where $\mathbf{1}_n$ is the indicator of $\{|b| \leq n\} \subset \mathbb R^{d+1}$.}
\end{equation}
We can  additionally mollify $b_n$ to obtain a $C^\infty$ smooth approximation of $b$ such that the Morrey norm of the approximating vector field does not exceed $(1+\varepsilon)\|b\|_{E_q}$ for any fixed $\varepsilon>0$, as is needed to turn a priori Sobolev regularity estimates for \eqref{eq1} into a posteriori estimates. However, a regularization of $b$ given by \eqref{b_n} will suffice (in particular, we will be able to apply It\^{o}'s formula to solutions of parabolic equations with drift $b_n$).

Put $b^{\frac{1}{p}}:=b|b|^{-1+\frac{1}{p}}$. 

Let $\mathcal E=\mathcal E_p:=\cup_{\varepsilon>0}e^{-\varepsilon |b|}L^p(\mathbb R^{d+1})$, a dense subspace of $L^p(\mathbb R^{d+1})$.

\medskip

\noindent{\sc Acknowledgements.} The author is grateful to Renming Song for some useful comments.

\bigskip

\section{Main results}

\label{par_sect}

\textbf{1.~}We first develop a Sobolev regularity theory of the inhomogeneous parabolic equation
\begin{equation}
\label{eq1}
(\lambda + \partial_t - \Delta + b(t,x) \cdot \nabla)u=f \quad \text{ on } \mathbb R^{d+1}.
\end{equation}
The next theorem is essential for the rest of the paper.

\begin{theorem}[Sobolev regularity theory]
\label{thm1}
Let
$
b=b_{\mathfrak s}+b_{\mathfrak b},
$
where 
\begin{equation}
\label{H}
\text{$|b_\mathfrak{s}| \in E_q$ for some $q>1$ close to $1$, and $|b_{\mathfrak b}| \in L^\infty(\mathbb R^{d+1})$} 
\end{equation}
(indices $\mathfrak s$ and $\mathfrak b$ stand for ``singular'' and ``bounded'', respectively).

The following are true:

{\rm(\textit{i})} For every $p \in ]1,\infty[$ there exist constants $c_{d,p,q}$ and $\lambda_{d,p,q}$ such that if
$$\|b_{\mathfrak s}\|_{E_q} < c_{d,p,q},$$ then,
for every $\lambda \geq \lambda_{d,p,q}$, solutions $u_n \in L^p(\mathbb R^{d+1})$ to the approximating parabolic equations
\begin{equation*}
(\lambda + \partial_t - \Delta + b_n \cdot \nabla)u_n=f, \quad f \in L^p(\mathbb R^{d+1})
\end{equation*}
converge in $\mathbb{W}^{1+\frac{1}{p},p}(\mathbb R^{d+1})$ to 
\begin{equation}
\label{u_repr}
u:=(\lambda+\partial_t-\Delta)^{-1}f - (\lambda+\partial_t-\Delta)^{-\frac{1}{2}-\frac{1}{2p}}Q_p (1+T_p)^{-1}R_p (\lambda+\partial_t-\Delta)^{-\frac{1}{2p'}} f,
\end{equation}
where the operators
$$R_p = R_p(b):=b^{\frac{1}{p}}\cdot \nabla (\lambda+\partial_t-\Delta)^{-\frac{1}{2}-\frac{1}{2p}},
$$
\begin{equation}
\label{Q_def}
Q_p = Q_p(b):=\bigl[(\lambda+\partial_t-\Delta)^{-\frac{1}{2p'}}|b|^{\frac{1}{p'}} \upharpoonright \mathcal E\bigr]^{\rm clos}_{p \rightarrow p}
\end{equation}
are bounded on $L^p(\mathbb R^{d+1})$,
and the operator 
$T_p:=R_p Q_p$ has norm 
\begin{equation}
\label{T_est}
\|T_p\|_{p \rightarrow p} <1.
\end{equation}

\smallskip

\rm{(}\textit{ii}\rm{)} If above $p>d+1$, then, by \eqref{u_repr} and by the parabolic Sobolev embedding, the convergence is uniform on $\mathbb R^{d+1}$ and $u \in C_\infty(\mathbb R^{d+1})$.
\end{theorem}

\begin{remarks}
1.~If $b$ is bounded, then $Q_p=(\lambda+\partial_t-\Delta)^{-\frac{1}{2p'}}|b|^{\frac{1}{p'}}$, and so the RHS of \eqref{u_repr} is simply the Duhamel series representation for the solution to \eqref{eq1} in $L^p(\mathbb R^{d+1})$ provided by the standard theory.

2.~The constraint $\|b_\mathfrak{s}\|_{E_q} < c_{d,p,q}$ is needed to ensure that $\|T_p\|_{p \rightarrow p} <1$, see Proposition \ref{prop2}. 

3.~If e.g.\,$b \in C_\infty(\mathbb R,L^d(\mathbb R^d))$ or $b=b(x)$ is in $L^d(\mathbb R^d)$, then one can represent $b=b_{\mathfrak s} + b_{\mathfrak b}$ with $\|b_{\mathfrak s}\|_{E_q}$ arbitrarily small (by defining $b_{\mathfrak b}$ to be a cutoff of $b$ such that the remaining part $b_{\mathfrak s}$ has sufficiently small $L^\infty(\mathbb R,L^d(\mathbb R^d))$ norm).
\end{remarks}

Given a general $b$ satisfying \eqref{H}, it is natural to ask, in what sense $u$ defined by \eqref{u_repr} solves the parabolic equation \eqref{eq1}?
Let us consider the case $p=2$ and $f \in L^2(\mathbb R^{d+1})$. 

\begin{definition*}
We say that a function $u \in \mathbb W^{\frac{3}{2},2}$ is a weak solution to \eqref{eq1} if the following identity is satisfied:
\begin{equation}
\label{weak_sol_def}
\begin{array}{r}
\langle (\lambda+\partial_t-\Delta)^{\frac{3}{4}}u, (\lambda + \partial_t-\Delta)^{\frac{3}{4}}\eta \rangle + \langle R_2(b) (\lambda+\partial_t-\Delta)^{\frac{3}{4}}u,Q_2^\ast(b)(\lambda+\partial_t-\Delta)^{\frac{3}{4}} \eta \rangle \\ [2mm]
=\langle f,(\lambda-\partial_t-\Delta)^{-\frac{1}{4}}(\lambda+\partial_t-\Delta)^{\frac{3}{4}}\eta \rangle
\end{array}
\end{equation}
for all $\eta \in C_c^\infty(\mathbb R^{d+1})$.
\end{definition*}

(Identity \eqref{weak_sol_def} is obtained by formally multiplying equation \eqref{eq1} by $(\lambda-\partial_t-\Delta)^{-\frac{1}{4}}(\lambda+\partial_t-\Delta)^{\frac{3}{4}}\eta$ and integrating over $\mathbb R^{d+1}$. Note that  $Q^\ast_2(b)=|b|^{\frac{1}{2}}(\lambda- \partial_t-\Delta)^{-\frac{1}{4}}$ is in $\mathcal B(L^2)$ by Proposition \ref{prop1}.)

Then $u$ defined by \eqref{u_repr}
is the unique in $\mathbb{W}^{\frac{3}{2},2}$ weak solution to \eqref{eq1}.
See Remark \ref{weak_sol_proof} for the proof.

\begin{remark}
\label{scale_shift_rem}
Thus, in order to have weak well-posedness of \eqref{eq1} for drifts satisfying \eqref{H} for $q>1$  we have to shift the standard scale of Hilbert spaces $\mathbb W^{1,2} \subset L^2(\mathbb R^{d+1}) \subset \mathbb W^{-1,2}$  to $\mathbb W^{\frac{3}{2},2} \subset \mathbb W^{\frac{1}{2},2} \subset \mathbb W^{-\frac{1}{2},2}$. If we were to consider instead of \eqref{eq1} a more general equation $(\lambda + \partial_t - \nabla \cdot a \cdot \nabla + b \cdot \nabla)u=f$ with a uniformly elliptic discontinuous matrix $a$, then the second order term would force us to work in the standard scale, and hence would require more restrictive assumptions on $b$: the form-boundedness, see the beginning of the paper (on the scale of Morrey spaces this will be \eqref{H} with $q>2$). 
\end{remark}

An analogous result can be obtained for $f \in L^p(\mathbb R^{d+1})$ with \eqref{weak_sol_def} modified according to \eqref{u_repr}.

\medskip

\textbf{2.~}Fix $T>0$. For given $n=1,2,\dots$ and $0 \leq r<T$, let $v_n \in C_b([r,T],C_\infty(\mathbb R^d))$ denote the solution to the Cauchy problem
\begin{equation}
\label{cauchy}
\left\{
\begin{array}{l}
(\lambda + \partial_t-\Delta + b_n(t,x) \cdot \nabla)v_n  =0 \quad (t,x) \in ]r,T] \times \mathbb R^d, \\[2mm]
 v_n(r,\cdot)=g(\cdot) \in C_\infty(\mathbb R^d),
\end{array}
\right.
\end{equation}
where $b_n$'s are defined by \eqref{b_n}. 
By a standard result, for every $n$, the operators $$U_n^{t,r}g:=v_n(t), \quad 0 \leq r \leq t \leq T$$ constitute a Feller evolution family on $C_\infty(\mathbb R^d)$.

Let $\delta_{s=r}$ denote the delta-function in the time variable $s$. Put
$$
(\lambda+\partial_t-\Delta)^{-1}\delta_{s=r}g(t,x):=\mathbf{1}_{t \geq r}e^{-\lambda (t-r)}(4\pi (t-r))^{-\frac{d}{2}}\int_{\mathbb R^d} e^{-\frac{|x-y|^2}{4(t-r)}} g(y)dy,
$$
$$
\nabla (\lambda+\partial_t-\Delta)^{-\frac{1}{2}-\frac{1}{2p'}}\delta_{s=r}g:=\mathbf{1}_{t \geq r}e^{-\lambda (t-r)} (t-r)^{-\frac{1}{2}+\frac{1}{2p'}} (4\pi (t-r))^{-\frac{d}{2}} \int_{\mathbb R^d} \nabla_x e^{-\frac{|x-y|^2}{4(t-r)}}g(y) dy.
$$

Recall $D_T=\{0 \leq r \leq t \leq T\}$.

\begin{theorem}[$C_\infty$ theory]
\label{thm2}
Under the assumptions of Theorem \ref{thm1}, let $\|b_{\mathfrak s}\|_{E_q} < c_{d,p,q}$ for a $p>d+1$. Then the following are true:

\smallskip

{\rm(\textit{i})} The limit
$$
U^{t,r}:=s\mbox{-}C_\infty(\mathbb R^d)\mbox{-}\lim_n U_n^{t,r}  \quad \text{uniformly in $(r,t) \in D_T$}
$$
exists and determines a Feller evolution family on $C_\infty(\mathbb R^d)$.

\smallskip

{\rm(\textit{ii})} For every initial function $g \in C_\infty(\mathbb R^d) \cap W^{1,p}(\mathbb R^d)$,  $v(t):=U^{t,r}g$, where $(r,t) \in D_T$, has representation
\begin{equation}
\label{v_repr}
v=(\lambda+\partial_t-\Delta)^{-1}\delta_{s=r}g - (\lambda+\partial_t-\Delta)^{-\frac{1}{2}-\frac{1}{2p}}Q_p (1+T_p)^{-1}G_p  S_p g,
\end{equation}
where
$
G_p=G_p(b):=b^{\frac{1}{p}} (\lambda+\partial_t-\Delta)^{-\frac{1}{2p}} \in \mathcal B(L^p(\mathbb R^{d+1}))
$
and
$
S_p g:=\nabla (\lambda+\partial_t-\Delta)^{-\frac{1}{2}-\frac{1}{2p'}}\delta_{s=r}g
$
satisfies 
$\|S_pg\|_{L^p(\mathbb R^{d+1})} \leq C_{p,d} \|\nabla g\|_{L^p(\mathbb R^d)}$.

\medskip

{\rm(\textit{iii})} As a consequence of \eqref{v_repr} and the parabolic Sobolev embedding, we obtain
$$
\sup_{(r,t) \in D_T, x \in \mathbb R^d}|v(t,x;r)| \leq C\|g\|_{W^{1,p}(\mathbb R^d)}.
$$
\end{theorem}

\medskip

\textbf{3.~}Recall that a probability measure $\mathbb P_x$, $x \in \mathbb R^d$ on $(C([0,T],\mathbb R^d),\mathcal B_t=\sigma(\omega_r \mid 0\leq r \leq t))$, where $\omega_t$ is the coordinate process, is said to be a martingale solution to SDE
\begin{equation}
\label{seq}
X_t=x-\int_0^t b(r,X_r)dr + \sqrt{2}B_t.
\end{equation} 
 if

1) $\mathbb P_x[\omega_0=x]=1$;

2) $\mathbb E_{x} \int_0^r |b(t,\omega_t)|dt<\infty$, $0<r \leq T$;

3) for every $f \in C_c^2(\mathbb R^d)$ the process
$$
r \mapsto f(\omega_r)-f(x) + \int_0^r (-\Delta f + b \cdot \nabla f)(t,\omega_t)dt
$$
is a $\mathcal B_r$-martingale under $\mathbb P_x$.

A martingale solution $\mathbb P_x$ of \eqref{seq} is said to be a weak solution if, upon completing $\mathcal B_t$ with respect to $\mathbb P_x$ (to, say, $\hat{\mathcal B}_t$), there exists a Brownian motion $B_t$ on $\big(C([0,T],\mathbb R^d),\hat{\mathcal B}_t,\mathbb P_x\big)$ such that
$$
\omega_r=x - \int_0^r b(t,\omega_t)dt + \sqrt{2}B_r, \quad r \geq 0 \quad \mathbb P_x\text{-a.s.}
$$

Put 
\begin{equation*}
P^{t,r}(b):=U^{T-t,T-r}(\tilde{b}), \qquad \tilde{b}(t,x)=b(T-t,x)
\end{equation*}
where $0 \leq t \leq r \leq T$.

\begin{theorem}
\label{thm3}
Under the assumptions of Theorem \ref{thm2} the following are true:

\smallskip

{\rm (\textit{i})} The backward Feller evolution family $\{P^{t,r}\}_{0 \leq t \leq r \leq T}$ is conservative, i.e.\,the density $P^{t,r}(x,\cdot)$ satisfies
$$ 
\langle P^{t,r}(x,\cdot)\rangle=1 \quad \text{ for all } x \in \mathbb R^d,
$$ 
and determines probability measures $\mathbb P_x$, $x \in \mathbb R^d$ on $(C([0,T],\mathbb R^d),\mathcal B_t)$, such that
$$
\mathbb E_x[f(\omega_r)]=P^{0,r}f(x), \quad 0 \leq r \leq T, \quad f \in C_\infty(\mathbb R^d).
$$

\smallskip

{\rm (\textit{ii})} For every $x \in \mathbb R^d$, the probability measure $\mathbb P_x$ is a weak solution to \eqref{seq}.

{\rm (\textit{iii})} For every $x \in \mathbb R^d$ and $\mathsf{f}$ satisfying \eqref{H}, given a $p>d+1$ as in Theorem \ref{thm2} (generally speaking, the larger $p$ is the smaller $\|b_{\mathfrak s}\|_{E_q}$ has to be), there exists constant $c$ such that for all $ h \in C_c(\mathbb R^{d+1})$ 
\begin{equation*}
\mathbb E_{x}\int_0^T |\mathsf{f}(r,\omega_r)h(r,\omega_r)|dr \leq c\|\mathbf{1}_{[0,T]}|\mathsf{f}|^{\frac{1}{p}} h\|_p
\end{equation*}
(in particular, one can take $\mathsf{f}=b$).

{\rm ($\textit{iii}$')} For every $x \in \mathbb R^d$, given a $\nu>\frac{d+2}{2}$, there exists a constant $c$ such that for all $ h \in C_c(\mathbb R^{d+1})$ the following Krylov-type bound is true:
\begin{equation}
\label{krylov_type}
\mathbb E_{x}\int_0^T |h(r,\omega_r)|dr \leq c\|\mathbf{1}_{[0,T]}h\|_\nu.
\end{equation}

\smallskip

{\rm (\textit{iv})} Any martingale solution $\mathbb Q_x$ to \eqref{seq} that satisfies for some $p>d+1$ as in Theorem \ref{thm2}
\begin{equation}
\label{kr_est1}
\mathbb E_{\mathbb Q_x}\int_0^T |b(r,\omega_r)h(r,\omega_r)|dr \leq c\|\mathbf{1}_{[0,T]}|b|^{\frac{1}{p}} h\|_p, \quad h \in C_c(\mathbb R^{d+1}), 
\end{equation}
coincides with $\mathbb P_x$.

\smallskip

{\rm ($\textit{iv}'$)} If additionally $|b| \in L_{\loc}^{\frac{d+2}{2}+\varepsilon}(\mathbb R^{d+1})$ for some $\varepsilon>0$ and $\|b_{\mathfrak s}\|_{E_q}$ is sufficiently small, then any martingale solution $\mathbb Q_x$ to \eqref{seq} that satisfies \eqref{krylov_type} for some $\nu>\frac{d+2}{2}$ sufficiently close to $\frac{d+2}{2}$ (depending on how small $\varepsilon$ is)
coincides with $\mathbb P_x$.

\end{theorem}

\begin{remark}
The uniqueness result ($\textit{iv}'$) is of the same type as e.g.\,in \cite{RZ}. 
\end{remark}

\subsection{Key bounds}

Here is the key bound used in the proofs of Theorems \ref{thm1} and \ref{thm2}.

\begin{proposition} 
\label{prop1}
Let $|b| \in E_q$ for some $q>1$ close to $1$.
Then for every $p \in ]1,\infty[$ there exists a constant $c_{p,q}$ such that
\begin{equation*}
\||b|^{\frac{1}{p}}(\pm\, \partial_t-\Delta)^{-\frac{1}{2p}}\|_{p \rightarrow p} \leq c_{p,q}\|b\|^{\frac{1}{p}}_{E_q}
\end{equation*}
\end{proposition}

In the time homogeneous case $b=b(x)$ the estimate on $\||b|^{\frac{1}{p}}(\lambda-\Delta)^{-\frac{1}{2p}}\|_{L^p(\mathbb R^d) \rightarrow L^p(\mathbb R^d)}$ in terms of the elliptic Morrey norm of $|b|$ is due to \cite[Theorem 7.3]{A}. Similar estimates in the parabolic case were obtained in \cite[proof of Prop.\,4.1]{Kr3}. There the proof is carried out for a different set of parameters than the one needed in this paper, so we included the details in Section \ref{prop1_proof_sect}.

As an immediate consequence of Proposition \ref{prop1} we obtain the following

\begin{proposition}
\label{prop2}
Let
$
b=b_\mathfrak{s}+b_{\mathfrak b},
$
satisfy \eqref{H}.
Then, for all $\lambda>0$,
\begin{equation}
\label{e1}
\|R_p\|_{p \rightarrow p} 
\leq C_{d,p,q}\|b_{\mathfrak{s}}\|^{\frac{1}{p}}_{E_q} + c \lambda^{-\frac{1}{2p}}\|b_{\mathfrak b}\|^{\frac{1}{p}}_{L^\infty(\mathbb R^{d+1})}
\end{equation}
\begin{equation}
\label{e2}
\|Q_p \|_{p \rightarrow p} \leq C'_{p,q}\|b_{\mathfrak s}\|^{\frac{1}{p'}}_{E_q}  + c'\lambda^{-\frac{1}{2p'}}\|b_{\mathfrak b}\|^{\frac{1}{p'}}_{L^\infty(\mathbb R^{d+1})}.
\end{equation}
\end{proposition}

The first estimate \eqref{e1} follows from Proposition \ref{prop1} 
using the boundedness of parabolic Riesz transforms (see \cite{G}).
The second estimate \eqref{e2} follows from Proposition \ref{prop1} by duality.

\subsection{Some remarks}

\begin{remark}
Theorems \ref{thm1}, \ref{thm2} are time inhomogeneous counterparts of the results in \cite{Ki}, Theorem \ref{thm3} is a time inhomogeneous counterpart of the result in \cite{KiS1}, where the authors treated vector fields $b:\mathbb R^d \rightarrow \mathbb R^d$ that can be more singular than the ones considered in this paper, namely, $|b| \in L^1_{\loc}(\mathbb R^d)$ and there exists $\delta>0$ such that, for some $\lambda>0$,
\begin{equation}
\label{wfbd}
\big\||b|^{\frac{1}{2}}\varphi\big\|_{L^2(\mathbb R^d)} \leq \delta \big\|(\lambda-\Delta)^{\frac{1}{4}}\varphi \big\|_{L^2(\mathbb R^d)}, \quad \forall\,\varphi \in C_c^\infty(\mathbb R^d).
\end{equation}
These vector fields are called weakly form-bounded. The fact that the time-homogeneous vector fields in $E_q$, $1<q \leq 2$ are weakly form-bounded follows from \cite[Theorem 7.3]{A}. A key difference between \cite{Ki} and the present paper is in the elliptic analogue of estimate \eqref{T_est}, i.e.
\begin{equation}
\label{T0_est}
\|\tilde{T}_p\|_{L^p(\mathbb R^d) \rightarrow L^p(\mathbb R^d)}<1, \quad \text{where }
\tilde{T}_p:=b^{\frac{1}{p}}\cdot \nabla(\mu-\Delta)^{-1} |b|^{\frac{1}{p'}},
\end{equation}
needed to ensure the convergence of the Neumann series.
This estimate is proved in \cite{Ki}
directly,
without splitting $\tilde{T}_p$ into a product of operators $b^{\frac{1}{p}}\cdot \nabla(\mu-\Delta)^{-\frac{1}{2}+\frac{1}{2p}}$ and $(\mu-\Delta)^{-\frac{1}{2p'}}|b|^{\frac{1}{p'}}$ as we do for $T_p$ in Theorem \ref{thm1}. In fact, the previous two operators are not even bounded on $L^p(\mathbb R^d)$ under the assumption \eqref{wfbd}; to have the $L^p(\mathbb R^d)$ boundedness one has to replace the exponents $-\frac{1}{2}+\frac{1}{2p}$, $-\frac{1}{2p'}$ by $-\frac{1}{2}+\frac{1}{2p}-\varepsilon$, $-\frac{1}{2p'}-\varepsilon$ for a $\varepsilon>0$ (this shows, by the way, that the difference between the elliptic Morrey class $M_q$ with $q>1$ and the larger class \eqref{wfbd} is quite significant). However, the proof of \eqref{T0_est} in \cite{Ki} requires inequalities for symmetric Markov generators, not valid for the parabolic operator $\partial_t-\Delta$ (although, it seems, one can address a parabolic counterpart of \eqref{wfbd} via an appropriate symmetrization of $\partial_t-\Delta$, which we plan to do in a subsequent paper).
\end{remark}

\begin{remark}
\label{a_rem}
Despite what was said in the introduction, \cite{KiM,Ki}  and the present paper are not directly comparable. Namely, the results in \cite{KiM,Ki} with $b \in \mathbf{F}_\delta$ can be extended more or less directly to the non-divergence form equation
\begin{equation}
\label{a_eq}
\big(-a(t,x)\cdot\nabla^2 + b(t,x) \cdot \nabla\big)u=f, \quad u(s,\cdot)=g(\cdot), \quad t>s,
\end{equation}
and the corresponding It\^{o}'s SDE having ``form-bounded diffusion coefficients''. Namely, the matrix $a$ is assumed to be bounded, uniformly elliptic and have 
\begin{equation}
\label{a_cond}
\nabla a_{ij} \in \mathbf{F}_{\delta_{ij}}, \quad 1 \leq i,j \leq d.
\end{equation}
See \cite{KiS2} where this scheme was carried out in the time-homogeneous case $b=b(x)$, $a=a(x)$.
(Let us note that since $\nabla a_{ij}(x)$ are form-bounded, they are in the elliptic Morrey class with $q=2$, see the introduction. Hence by Poincar\'{e}'s inequality such $a_{ij}$ belong to the ${\rm VMO}$ class, see details in \cite[Sect.3]{Kr0}.) 
However, a similar extension of the results of the present paper for \eqref{sde1} with $b \in E_{q}$, $1<q\leq 2$ requires more regular $a$, see e.g.\,Remark \ref{scale_shift_rem}.
\end{remark}

\bigskip

\section{Some corollaries of Theorem \ref{thm1} and Theorem \ref{thm2}(\textit{ii})}

\label{cor_sect}

The following results will be needed in the proof of Theorem \ref{thm3}. 

\begin{corollary}
\label{cor1}
Let vector fields $b$, $\mathsf{f}$ satisfy \eqref{H}. Let $b_n$ be given by \eqref{b_n} and let us define $\mathsf{f}_n$ similarly. Then, under the assumptions on $\|b_{\mathfrak s}\|_{E_q}$ of Theorem \ref{thm1}, 
for every $\lambda \geq \lambda_{d,p,q}$ solutions $u_n \in L^p(\mathbb R^{d+1})$ to the approximating parabolic equations
$$
(\lambda + \partial_t - \Delta + b_n \cdot \nabla)u_n=|\mathsf{f}_n|h, \quad h \in C_c(\mathbb R^{d+1})
$$
converge in $\mathbb{W}^{1+\frac{1}{p},p}(\mathbb R^{d+1})$ to 
\begin{equation*}
u:=(\lambda+\partial_t-\Delta)^{-1}|\mathsf{f}|h - (\lambda+\partial_t-\Delta)^{-\frac{1}{2}-\frac{1}{2p}}Q_p(b) (1+T_p(b))^{-1}R_p(b) Q_p(\mathsf{f})|\mathsf{f}|^{\frac{1}{p}}h,
\end{equation*}
In particular, 
if $p>d+1$, then the convergence is uniform on $\mathbb R^{d+1}$ and
$$
\sup_{\mathbb R^{d+1}}|u| \leq C \|\mathsf{f}^{\frac{1}{p}}h\|_p.
$$

\end{corollary}

The reason we include $h$ in Corollary \ref{cor1} is because in general $|\mathsf{f}|$ is only in $L^1_\loc(\mathbb R^{d+1})$, not in $L^1(\mathbb R^{d+1})$. In the proof of Theorem \ref{thm3} we will be applying Corollary \ref{cor1} and Corollaries \ref{cor2}, \ref{cor3} below to $\mathsf{f}_n$ equal to either $b_n$ or (with some abuse of notation) to $b_n-b_k$.

\begin{corollary}
\label{cor2}
Under the assumptions and notation of Corollary \ref{cor1}, if $p>d+1$, then, for every $\lambda \geq \lambda_{d,p,q}$, solutions $v_n \in C_b([r,T],C_\infty(\mathbb R^d))$ to the approximating inhomogeneous Cauchy problems
\begin{align*}
(\lambda + \partial_t - \Delta + b_n \cdot \nabla)v_n & =\mathbf{1}_{[r,T]}|\mathsf{f}_n|h \quad \text{ on } ]r,T] \times \mathbb R^d, \\
 v_n(r,\cdot)&=g(\cdot) \in C_\infty(\mathbb R^d) \cap W^{1,p}(\mathbb R^d),
\end{align*}
where $0 \leq r<T$,
converge uniformly on $D_T \times \mathbb R^d$ to 
\begin{align*}
v& :=(\lambda+\partial_t-\Delta)^{-1}\bigl[\mathbf{1}|\mathsf{f}|h + \delta_{s=r}g\bigr] \\
&- (\lambda+\partial_t-\Delta)^{-\frac{1}{2}-\frac{1}{2p}}Q_p(b) (1+T_p(b))^{-1} \bigl[R_p(b) Q_p(\mathbf{1}\mathsf{f})|\mathbf{1}\mathsf{f}|^{\frac{1}{p}}h + G_p(b)  S_p g\bigr], \quad \mathbf{1}:=\mathbf{1}_{[r,T]},
\end{align*}
$$
\sup_{(r,t) \in D_T, x \in \mathbb R^d}|v(t,x;r)| \leq C_1 \|\mathbf{1}\mathsf{f}^{\frac{1}{p}}h\|_p + C_2\|g\|_{W^{1,p}(\mathbb R^d)}.
$$
\end{corollary}

We also have the following weighted variant of Corollary \ref{cor2}, which we will record for bounded vector fields $b_n$ and $\mathsf{f}_n$ and $\lambda=0$, as will be needed in the proof of Theorem \ref{thm3} in Section \ref{proof_thm3_sect}.
 
Let $q>1$ (close to $1$) be from the hypothesis on $b$ in Theorem \ref{thm1}.
Set
$$\rho(x):=(1+l |x|^2)^{-\nu},   \quad x \in \mathbb R^d,$$
where $\nu > \frac{d}{2p}+\frac{1}{pq'}$ is fixed (so that $\rho \in L^p(\mathbb R^d)$ and \eqref{rho_fin} below holds) and $l>0$ is to be chosen. We have
\begin{equation}
\label{eta_two_est}
\tag{$\star$}
|\nabla \rho| \leq \nu \sqrt{l}\rho =: c_1 \sqrt{l}\rho, \quad |\Delta \rho| \leq 2\nu (2\nu + d+2 ) l \rho =: c_2 l \rho.
\end{equation}

\begin{corollary}
\label{cor3}
Under the assumptions  of Corollary \ref{cor1}, if $p>d+1$, then, provided that the constant $l$ in the definition of  $\rho$ is chosen sufficiently small, solutions $v_n \in C_b([r,T],C_\infty(\mathbb R^d))$ to the approximating inhomogeneous Cauchy problems
\begin{align*}
(\partial_t - \Delta + b_n \cdot \nabla)v_n & =\pm\mathbf{1}_{[r,T]}|\mathsf{f}_n| \quad \text{ on } [r,T] \times \mathbb R^d, \\
 v_n(r,\cdot)&=g \in C_\infty(\mathbb R^d) \cap W^{1,p}(\mathbb R^d),
\end{align*}
where $0 \leq r<T$,
satisfy, for all $t \in ]r,T]$,
\begin{equation}
\label{weight_bd0}
\sup_{[r,t] \times \mathbb R^d}|\rho v_n| \leq C_1 \|\rho\mathbf{1}_{[r,t]}|\mathsf{f}_n|^{\frac{1}{p}}\|_p + C_2\|\rho g\|_{W^{1,p}(\mathbb R^d)}
\end{equation}
and, putting $\rho_y(x):=\rho(x-y)$,
\begin{align}
\label{weight_bd}
\sup_{[r,t] \times \mathbb R^d}|v_n| & \leq \sup_{y \in \mathbb Z^d}(C_1\|\rho_y\mathbf{1}_{[r,t]}|\mathsf{f}_n|^{\frac{1}{p}}\|_p + C_2\|\rho_y g\|_{W^{1,p}(\mathbb R^d)}) \\
& \leq \tilde{C}_1(t-r)^\gamma\big(\|\mathsf{f}_{\mathfrak s}\|_{E_q}^{\frac{1}{p}} + \|\mathsf{f}_{\mathfrak b}\|_{\infty}^{\frac{1}{p}}\big) + \tilde{C}_2 \|g\|_{W^{1,p}} \label{weight_bd2}
\end{align}
with constants $C_1$, $C_2$, $\tilde{C}_1$, $\tilde{C}_2$ and $\gamma>0$ independent of $n$ and $t$.
\end{corollary}

We prove Corollary \ref{cor3} in Section \ref{cor3_proof_sect}.

\bigskip

\section{Proof of Proposition \ref{prop1}}

\label{prop1_proof_sect}

It suffices to carry out the proof for $b_n$ defined by \eqref{b_n} and then use the Dominated Convergence Theorem. So, without loss of generality, everywhere below $b$ is bounded. Below we follow \cite{Kr3}.

Set
$$
M_\beta h(t,x):=\sup_{\rho>0}\rho^\beta \frac{1}{|C_\rho|}\int_{C_\rho(t,x)} |h| dz, \quad 0 \leq \beta \leq d-2,
$$
and define the maximal function $$M h(t,x):=\sup_{\rho>0} \frac{1}{|C_\rho|}\int_{C_\rho(t,x)} |h| dz.$$ 
Also, define 
$$
\hat{M} h(t,x):=\sup_{(t,x) \in C} \frac{1}{|C|}\int_{C} |h| dz,
$$
where the supremum is taken over all cylinders $C=C_\rho(z) \ni (t,x)$, $z \in \mathbb R^{d+1}$, $\rho>0$.

Put $P_\alpha:=(-\partial_t-\Delta)^{-\frac{\alpha}{2}}$.  We will need

\begin{lemma}[see {\cite[Lemma 2.2]{Kr3}}] 
\label{lem_m}
For every $\beta \in ]0,d+2]$, $0<\alpha<\beta$ there exists $C>0$ such that, for all $f \geq 0$,
$$
P_\alpha f \leq C (M_\beta f)^{\frac{\alpha}{\beta}}(Mf)^{1-\frac{\alpha}{\beta}}.
$$
\end{lemma}

Let us prove the first inequality:
\begin{equation}
\label{req_ineq}
|\langle |b|(P_{\frac{1}{p}}f)^p\rangle| \leq c_{p,q}^p\|b\|_{E_q}\|f\|_p^p, \quad f \in L^p(\mathbb R^{d+1})
\end{equation}
(the proof of the second one is similar). 
Put $u:=P_{\frac{1}{p}}f$. Then we estimate the LHS of \eqref{req_ineq} as
\begin{align}
|\langle |b|(P_{\frac{1}{p}}f)^p\rangle| &= |\langle |b||u|^{p-1},P_{\frac{1}{p}}f\rangle| \notag \\
& \leq |\langle P^\ast_{\frac{1}{p}}(|b||u|^{p-1}),f\rangle| \leq \|P^\ast_{\frac{1}{p}}(|b||u|^{p-1})\|_{p'}\|f\|_p. \label{ineq2}
\end{align}
Here $P^\ast_{\alpha}=(\partial_t-\Delta)^{-\frac{\alpha}{2}}$ is the adjoint of the operator $P_{\alpha}$ in $L^2(\mathbb R^{d+1})$. To obtain \eqref{req_ineq}, we need to bound the coefficient $\|P^\ast_{\frac{1}{p}}(|b|u^{p-1})\|_{p'}$. 
To this end, we estimate pointwise
\begin{align*}
P_{\frac{1}{p}}^\ast(|b||u|^{p-1}) = P_{\frac{1}{p}}^\ast(|b|^{\frac{1}{p}+\gamma} |b|^{\frac{1}{p'}-\gamma}|u|^{p-1}) \leq P_{\frac{1}{p}}^\ast(|b|^{1+\gamma p})^{\frac{1}{p}} (P_{\frac{1}{p}}^\ast(|b|^{1-\gamma p'}|u|^p))^{\frac{1}{p'}}
\end{align*} 
for a small $\gamma>0$ such that $1+\gamma p < q_0$ for some fixed $q_0<q$. 
Hence 
\begin{align}
\|P^\ast_{\frac{1}{p}}(|b||u|^{p-1})\|_{p'}^{p'} \leq \langle |b|^{1-\gamma p'}|u|^p, P_{\frac{1}{p}}[P^\ast_{\frac{1}{p}}(|b|^{1+\gamma p})]^{\frac{1}{p-1}}\rangle. \label{PP}
\end{align}

1) By Lemma \ref{lem_m} with $\alpha=\frac{1}{p}$, $\beta=1+\gamma p$ (or rather its straightforward variant for $P_\alpha^\ast$),
\begin{align*}
P^\ast_{\frac{1}{p}}(|b|^{1+\gamma p}) 
& \leq C\|b\|_{E_{1+\gamma p}}^{\frac{1}{p}}(\hat{M}|b|^{1+\gamma p})^{1-\frac{1}{p}\frac{1}{1+\gamma p}} \\
& \leq C\|b\|_{E_{q_0}}^{\frac{1}{p}}(\hat{M}|b|^{1+\gamma p})^{1-\frac{1}{p}\frac{1}{1+\gamma p}}
\end{align*}
At this point let us assume that
\begin{equation}
\label{a1}
\hat{M}|b|^{1+\gamma p} \leq C_0 |b|^{1+\gamma p},
\end{equation}
but will get rid of this assumption later. Then
$$
P^\ast_{\frac{1}{p}}(|b|^{1+\gamma p}) \leq C_2\|b\|_{E_{q_0}}^{\frac{1}{p}}|b|^{1+\gamma p - \frac{1}{p}}.
$$

2) After applying the last estimate in \eqref{PP}, it remains to estimate
$
P_{\frac{1}{p}} (|b|^{\frac{1}{p}+\gamma p'}).
$
Selecting $\gamma$ even smaller, if needed, one may assume that $\frac{1}{p}+\gamma p'<q_0$.
By Lemma \ref{lem_m},
\begin{align*}
P_{\frac{1}{p}} (|b|^{\frac{1}{p}+\gamma p'}) & \leq C (M_{\frac{1}{p}+\gamma p'}|b|^{\frac{1}{p}+\gamma p'})^{\frac{1}{p}\frac{1}{\frac{1}{p}+\gamma p'}}(M|b|^{\frac{1}{p}+\gamma p'})^{1-\frac{1}{p}\frac{1}{\frac{1}{p}+\gamma p'}} \\
& \leq C\|b\|_{E_{q_0}}^{\frac{1}{p}} (\hat{M}|b|^{\frac{1}{p}+\gamma p'})^{1-\frac{1}{p}\frac{1}{\frac{1}{p}+\gamma p'}}
\end{align*}
In addition to \eqref{a1}, let us temporarily assume that
\begin{equation}
\label{a2}
\hat{M}|b|^{\frac{1}{p}+\gamma p'} \leq C_0|b|^{\frac{1}{p}+\gamma p'}.
\end{equation}
Then
$$
P_{\frac{1}{p}} (|b|^{\frac{1}{p}+\gamma p'})  \leq C_2 \|b\|_{E_{q_0}}^{\frac{1}{p}} |b|^{\gamma p'}.
$$

\medskip

3) Applying the results from 1), 2)  in \eqref{PP}, we obtain
$$
\|P^\ast_{\frac{1}{p}}(|b||u|^{p-1})\|_{p'} \leq C_3 \|b\|_{E_{q_0}}^{\frac{1}{p}} \langle |b||u|^p\rangle^{\frac{1}{p'}}.
$$
Therefore, \eqref{ineq2} yields
$$
\langle |b||u|^{p}\rangle \leq C_3 \|b\|_{E_{q_0}}^{\frac{1}{p}} \langle |b||u|^p\rangle^{\frac{1}{p'}}\|f\|_p
$$
so 
\begin{equation}
\label{ineq4}
\langle |b||u|^{p}\rangle^{\frac{1}{p}} \leq C_3\|b\|_{E_{q_0}}^{\frac{1}{p}} \|f\|_p.
\end{equation}

\medskip

4) Now one gets rid of the assumptions \eqref{a1} and \eqref{a2} at expense of replacing $\|b\|_{E_{q_0}}$ in \eqref{ineq4} by $\|b\|_{E_q}$, where, recall, $q_0<q$. This, in turn, will give \eqref{req_ineq}.
Fix $q_0<q_1<q$ and define $\tilde{b}:=(\hat{M} |b|^{q_1})^{\frac{1}{q_1}}$. Then $\tilde{b} \geq |b|$ and $\tilde{b}$ satisfies 
\begin{equation}
\label{a3}
\hat{M}\tilde{b}^{q_0} \leq C_0\tilde{b}^{q_0}
\end{equation}
(see \cite[p.158]{GR}).
Since $1+\gamma p < q_0$, $\frac{1}{p}+\gamma p'<q_0$, both inequalities \eqref{a1} and \eqref{a2} for $\tilde{b}$ follow from \eqref{a3}, and so we have
$$
\langle |b||u|^{p}\rangle^{\frac{1}{p}} \leq C_3 \|\tilde{b}\|_{E_{q_0}}^{\frac{1}{p}} \|f\|_p.
$$
It remains to apply inequality $\|\tilde{b}\|_{E_{q_0}} \leq C \|b\|_{E_{q}}$, which was established in \cite[proof of Prop.\,4.1]{Kr3}. \hfill \qed

\bigskip

\section{Proof of Theorem \ref{thm1}}

The fact that $Q_p$, $R_p$ are bounded on $L^p(\mathbb R^{d+1})$, and the operator $T_p=R_pQ_p$ has norm $\|T_p\|_{p \rightarrow p}<1$ provided that $\|b_{\mathfrak s}\|_{E_q}$ is sufficiently small and $\lambda$ is sufficiently large, is an immediate consequence of Proposition \ref{prop2}. Thus, 
$$
\Theta_p(b):=(\lambda+\partial_t-\Delta)^{-1} - (\lambda+\partial_t-\Delta)^{-\frac{1}{2}-\frac{1}{2p}}Q_p (1+T_p)^{-1}R_p (\lambda+\partial_t-\Delta)^{-\frac{1}{2p'}}
$$
is in $\mathcal B(L^p)$.

Recall: $b_n:=\mathbf{1}_n b$, where $\mathbf{1}_n$ is the indicator of $\{(t,x) \in \mathbb R^{d+1} \mid |b(t,x)| \leq n\}$. If $u_n \in L^p(\mathbb R^{d+1})$ denotes the solution to
$
(\partial_t-\Delta + b_n \cdot \nabla)u_n=f,
$
which exists by the classical theory, then $$u_n=\Theta_p(b_n)f,$$ where $\Theta_p(b_n)f$ coincides with the Duhamel series representation for $u_n$.

Next, let us note that 
\begin{equation}
\label{R_Q_conv}
R_p(b_n) \rightarrow R_p(b), \;\; Q_p(b_n) \rightarrow Q_p(b) \quad \text{ strongly in $L^p(\mathbb R^{d+1})$}
\end{equation}
 as follows from \eqref{e1}, \eqref{e2} and the Dominated Convergence Theorem. Hence 
\begin{equation}
\label{u_conv}
u_n:=\Theta_p(b_n) \rightarrow u=\Theta_p(b) \text{ in  }\mathbb W^{1+\frac{1}{p},p}(\mathbb R^{d+1}),
\end{equation}
as needed.
\hfill \qed

\begin{remark}
\label{weak_sol_proof}
In the comment after Theorem \ref{thm1} we promised to prove the existence and uniqueness of weak solution to \eqref{eq1}. The argument is standard and goes as follows.
In the proof of Theorem \ref{thm1} above take $p=2$, so $u \in \mathbb W^{\frac{3}{2},2}$. Multiplying 
$(\lambda + \partial_t - \Delta + b_n \cdot \nabla)u_n=f$, $n=1,2,\dots$, by $\varphi=(\lambda - \partial_t-\Delta)^{-\frac{1}{4}}(\lambda + \partial_t-\Delta)^{\frac{3}{4}}\eta$, $\eta \in C_c^\infty(\mathbb R^{d+1})$ and integrating over $\mathbb R^{d+1}$, we have 
\begin{align*}
 \langle (\lambda+\partial_t-\Delta)^{\frac{3}{4}}u_n, (\lambda + \partial_t-\Delta)^{\frac{3}{4}}\eta\rangle
+ & \langle R_2(b_n)(\lambda+\partial_t-\Delta)^{\frac{3}{4}} u_n,Q_2^\ast(b_n)(\lambda + \partial_t-\Delta)^{\frac{3}{4}}\eta \rangle \\
&=\langle f,(\lambda - \partial_t-\Delta)^{-\frac{1}{4}}(\lambda + \partial_t-\Delta)^{\frac{3}{4}}\eta\rangle,
\end{align*}
where, recall, $Q^\ast_2(b)=|b|^{\frac{1}{2}}(\lambda- \partial_t-\Delta)^{-\frac{1}{4}} \in \mathcal B(L^2)$.
In view of \eqref{u_conv},
$$
\langle (\lambda+\partial_t-\Delta)^{\frac{3}{4}}u_n, (\lambda + \partial_t-\Delta)^{\frac{3}{4}}\eta\rangle \rightarrow \langle (\lambda+\partial_t-\Delta)^{\frac{3}{4}}u, (\lambda + \partial_t-\Delta)^{\frac{3}{4}}\eta\rangle \quad (n \rightarrow \infty).
$$
Next,
\begin{align*}
\langle R_2(b_n)(\lambda+\partial_t-\Delta)^{\frac{3}{4}} u_n & ,Q_2^\ast(b_n)(\lambda + \partial_t-\Delta)^{\frac{3}{4}}\eta\rangle \\
& = \langle R_2(b_n)(\lambda+\partial_t-\Delta)^{\frac{3}{4}} (u_n-u),Q_2^\ast(b_n)(\lambda + \partial_t-\Delta)^{\frac{3}{4}}\eta\rangle \\  
& + \langle R_2(b_n)(\lambda+\partial_t-\Delta)^{\frac{3}{4}} u,(Q_2^\ast(b_n)-Q_2^\ast(b))(\lambda + \partial_t-\Delta)^{\frac{3}{4}}\eta\rangle \\
& + \langle R_2(b_n)(\lambda+\partial_t-\Delta)^{\frac{3}{4}} u,Q_2^\ast(b)(\lambda + \partial_t-\Delta)^{\frac{3}{4}}\eta\rangle.
\end{align*}
By \eqref{u_conv}, $Q_2^\ast(b_n) \rightarrow Q_2^\ast(b)$ strongly in $L^2(\mathbb R^{d+1})$ (by the same argument as in \eqref{R_Q_conv}) we get that the first two terms in the RHS tend to $0$ as $n \rightarrow \infty$. By \eqref{R_Q_conv}, the last term tends to $\langle R_2(b)(\lambda+\partial_t-\Delta)^{\frac{3}{4}} u,Q_2^\ast(b)(\lambda + \partial_t-\Delta)^{\frac{3}{4}}\eta\rangle$. Hence $u$ is a weak solution to \eqref{eq1} in the sense of definition \eqref{weak_sol_def}. 

Let $v\in \mathbb W^{\frac{3}{2},2}$ be another weak solution. Put
$$
\tau[v,\eta]:=\langle (\lambda+\partial_t-\Delta)^{\frac{3}{4}}v, (\lambda+\partial_t-\Delta)^{\frac{3}{4}}\eta \rangle + \langle R_2(b)(\lambda+\partial_t-\Delta)^{\frac{3}{4}}v,Q^\ast_2(b)(\lambda+\partial_t-\Delta)^{\frac{3}{4}}\eta\rangle,
$$
where $\eta \in C_c^\infty(\mathbb R^{d+1})$. We have
$$
|\langle R_2(b)(\lambda+\partial_t-\Delta)^{\frac{3}{4}}v,Q^\ast_2(b)(\lambda+\partial_t-\Delta)^{\frac{3}{4}}\eta\rangle| \leq c\|v\|_{\mathbb W^{\frac{3}{2},2}}\|\eta\|_{\mathbb W^{\frac{3}{2},2}}
$$
where $c<1$ by our assumption on $b$. We extend $\tau[v,\eta]$ to $\eta \in \mathbb W^{\frac{3}{2},2}$ by continuity. Now we have $\tau[v-u,\eta]=0$, where $u$ is the weak solution constructed above, so it suffices to choose $\eta=v-u$ to arrive at $
0=\tau[v-u,v-u] \geq (1-c)\|v\|_{\mathbb W^{\frac{3}{2},2}}^2,
$ hence $v=u$.
\end{remark}

\bigskip

\section{Proof of Theorem \ref{thm2}}

Let us define $U^{t,r}g:=v(t)$ ($t \geq r$), $g \in C_\infty(\mathbb R^d) \cap W^{1,p}(\mathbb R^d)$ where $v(t)$ is given by \eqref{v_repr}. 
Since $U^{t,r}_n$ are $L^\infty$ contractions, it suffices to prove (we consider convergence in $(r,t) \in D_T$)
\begin{equation}
\label{v_conv}
U=C_b(D_T,C_\infty(\mathbb R^d))\mbox{-}\lim_{n}U_n g,
\end{equation}
and then extend operators $U^{t,r}$ by continuity to $g \in C_\infty(\mathbb R^d)$. The reproduction property of $U^{t,r}$ and the preservation of positivity will follow from the corresponding properties of $U^{t,r}_n$.

Proof of \eqref{v_conv}. Put $v_n:=U^{t,r}_n g$.
We have
$$
v_n=(\lambda+\partial_t-\Delta)^{-1}\delta_{s=r}g - (\lambda+\partial_t-\Delta)^{-\frac{1}{2}-\frac{1}{2p}}Q_p(b_n) (1+T_p(b_n))^{-1}G_p(b_n)  S_p g.
$$
This is the usual Duhamel series representation for $v_n$.
We know from the proof of Theorem \ref{thm1} that operators $Q_p(b_n)$, $T_p(b_n)$, $G_p(b_n)$ are bounded on $L^p(\mathbb R^{d+1})$ with operator norms independent of $n$. In turn, operator $S_p$
satisfies
$$
\|S_pg\|_{L^p(\mathbb R^{d+1})} \leq C_{p,d} \|\nabla g\|_{L^p(\mathbb R^d)}.
$$
(Indeed, taking for brevity $r=0$, we have by definition
$$
S_p g(t,x)=\mathbf{1}_{t \geq 0}e^{-\lambda t} t^{-\frac{1}{2}+\frac{1}{2p'}} (4\pi t)^{-\frac{d}{2}} \int_{\mathbb R^d} \nabla_x e^{-\frac{|x-y|^2}{4t}}g(y) dy.
$$
Hence
$$
\|S_p g\|_{L^p(\mathbb R^{d+1})}^p=\int_{\mathbb R} \mathbf{1}_{t \geq 0}\|S_p(t)\|_{L^p(\mathbb R^{d})}^p dt \leq \int_0^\infty e^{-\lambda p t} t^{(-\frac{1}{2}+\frac{1}{2p'})p} dt \|\nabla g\|_{L^p(\mathbb R^d)}^p,
$$
where $(-\frac{1}{2}+\frac{1}{2p'})p=-\frac{1}{2}$ $(>-1$ so the integral in time converges).)

Clearly, $(\lambda+\partial_t-\Delta)^{-1}\delta_{s=r}g \in C_b([r,\infty[,C_\infty(\mathbb R^d))$.
Thus, to prove \eqref{v_conv}, it remains to note that $Q_p(b_n) \rightarrow Q_p(b)$, $T_p(b_n) \rightarrow T_p(b)$ and $G_p(b_n) \rightarrow G_p(b)$ strongly in $L^{p}(\mathbb R^{d+1})$ (see proof of Theorem \ref{thm1}), so that by the parabolic Sobolev embedding, since $p>d+1$, 
\begin{align*}
&(\lambda+\partial_t-\Delta)^{-\frac{1}{2}-\frac{1}{2p}}Q_p(b_n) (1+T_p(b_n))^{-1}G_p(b_n)  S_p g \\
&\rightarrow (\lambda+\partial_t-\Delta)^{-\frac{1}{2}-\frac{1}{2p}}Q_p(b) (1+T_p(b))^{-1}G_p(b)  S_p g \quad \text{in $C_\infty(\mathbb R^{d+1})$ as $n \rightarrow \infty$.} 
\end{align*}
The required convergence \eqref{v_conv} follows. \hfill \qed

\bigskip

\section{Proof of Corollary \ref{cor3} (weighted bounds)}
\label{cor3_proof_sect}

It will be convenient to carry out the proof for solutions $v_n$ to  
\begin{align*}
(\lambda + \partial_t - \Delta + b_n \cdot \nabla)v_n & =\pm\mathbf{1}_{[r,T]}|\mathsf{f}_n|, \quad \text{ on } [r,T] \times \mathbb R^d, \\
 v_n(r,\cdot)&=g \in C_\infty(\mathbb R^d) \cap W^{1,p}(\mathbb R^d)
\end{align*}
where $\lambda \geq \lambda_{d,p,q}>0$. Since we are working on a fixed finite time interval $[r,T]$, this change will amount to multiplying solution $v_n$ by a bounded function $e^{\lambda(t-r)}$ which, clearly, does not affect the sought estimates.   

Put for brevity $v=v_n$. We will carry out the proof on the maximal interval $[r,t]$, i.e.\,for $t=T$. Also, without loss of generality, the RHS of the equation is $\mathbf{1}_{[r,T]}|\mathsf{f}_n|$ and the initial function $g \geq 0$, so $v \geq 0$.  We have for the weight $\rho$ defined before the corollary,
\begin{align*}
& \lambda \rho v +\partial_t (\rho v) - \Delta (\rho v) + b_n \cdot \nabla (\rho v) \\
& = \rho(\lambda v + \partial_t v - \Delta v +  b_n \cdot \nabla v) - 2 \nabla \rho \cdot \nabla v + (-\Delta \rho)v + b_n v \cdot \nabla \rho \\
& (\text{$v$ solves the parabolic equation above}) \\
& = \rho \mathbf{1}_{[r,T]}|\mathsf{f}_n| +K, \quad K:= - 2 \nabla \rho \cdot \nabla v + (-\Delta \rho)v + b_n v \cdot \nabla \rho.
\end{align*}
Let us rewrite the term $K$ as follows:
$$
K=-2\biggl(\frac{\nabla \rho}{\rho}\cdot \nabla (\rho v) - \frac{(\nabla \rho)^2}{\rho} v \biggr) + (-\Delta \rho)v + b_n v \cdot \nabla \rho.
$$
Hence 
\begin{equation}
\label{eq7}
\lambda \rho v +\partial_t (\rho v) - \Delta (\rho v) + \tilde{b}_n \cdot \nabla (\rho v) = \rho \mathbf{1}_{[r,T]}|\mathsf{f}_n| + \tilde{K},
\end{equation}
where $\tilde{b}_n:=b_n+2\frac{\nabla \rho}{\rho}$ and
$$
\tilde{K}=2 \frac{(\nabla \rho)^2}{\rho} v + (-\Delta \rho)v + b_n v \cdot \nabla \rho.
$$
Note that by \eqref{eta_two_est} the term $2\frac{\nabla \rho}{\rho}$ in $\tilde{b}_n$ is a bounded vector field whose $\sup$ norm can be made as small as needed by selecting $l$ in the definition of $\rho$ sufficiently small; we will be selecting $l$ sufficiently small.

By \eqref{eta_two_est}, the first two terms in $\tilde{K}$ are smaller than $(c_1^2+c_2) l\rho v $, so they can be handled by selecting $\lambda \geq \lambda_{d,p,q} + (c_1^2+c_2) l$. Corollary \ref{cor2}  applied to \eqref{eq7} gives us
$$
\sup_{[r,T] \times \mathbb R^d}|\rho v| \leq C_1 \|\rho\mathbf{1}_{[r,T]}|\mathsf{f}_n|^{\frac{1}{p}}\|_p + C_1'\sqrt{l}\||b_n|^{\frac{1}{p}}\rho v \|_p + C_2\|\rho g\|_{W^{1,p}(\mathbb R^d)},
$$
where, when we were treating the last term in $\tilde{K}$, we used again \eqref{eta_two_est}.
However, now we have to deal with $\||b_n|^{\frac{1}{p}}\rho v \|_p$. So, we will have to use instead a finer consequence of the solution representation of Corollary \ref{cor2}:
\begin{align}
\label{v_est}
\frac{1}{2}\sup_{[r,T] \times \mathbb R^d}|\rho v| &  + \frac{C_0}{2}\|(\lambda+\partial_t-\Delta)^{\frac{1}{2p}}\rho v\|_p \notag \\
& \leq C_1 \|\rho\mathbf{1}_{[r,T]}|\mathsf{f}_n|^{\frac{1}{p}}\|_p + C_1'\sqrt{l}\||b_n|^{\frac{1}{p}}\rho v\|_p + C_2\|\rho g\|_{W^{1,p}(\mathbb R^d)}
\end{align}
for appropriate constant $C_0>0$. (Here we only need to justify that $(\lambda+\partial_t-\Delta)^{-\frac{1}{2}-\frac{1}{2p'}}\delta_{s=r}\rho g$ is in $L^p(\mathbb R^{d+1})$ and its norm is bounded by $\|\rho g\|_{L^p(\mathbb R^d)}$. The proof of this fact repeats the argument for operator $S_p$ from the proof of Theorem \ref{thm2}.) We have
\begin{align*}
\||b_n|^{\frac{1}{p}}\rho v\|_p & \leq \||b_n|^{\frac{1}{p}}(\lambda+\partial_t-\Delta)^{-\frac{1}{2p}}\|_{p \rightarrow p}\|(\lambda+\partial_t-\Delta)^{\frac{1}{2p}}\rho v\|_p \\
& (\text{we are applying Proposition \ref{prop2}}) \\
& \leq C_\lambda \|(\lambda+\partial_t-\Delta)^{\frac{1}{2p}}\rho v\|_p.
\end{align*}
Applying the last bound in \eqref{v_est} with $l$ chosen sufficiently small so that $\frac{C_0}{2}-C_1'C_\lambda \sqrt{l}>0$, we obtain \eqref{weight_bd0}. 

Assertion \eqref{weight_bd} follows from \eqref{weight_bd0} using translations.

Armed with \eqref{weight_bd}, we now prove \eqref{weight_bd2}. We obtain from \eqref{eta_two_est} that $\sup_{y \in \mathbb Z^d}\|\rho_y g\|_{W^{1,p}} \leq c_0\|g\|_{W^{1,p}}$ for appropriate $c_0>0$.
It remains to show that if $|\mathsf{f}| \in E_q$, $q>1$, then 
\begin{equation}
\label{rho_fin}
\sup_{y \in \mathbb Z^d}\|\rho_y\mathbf{1}_{[r,T]}|\mathsf{f}|^{\frac{1}{p}}\|_p \leq c(T-r)^\gamma\|\mathsf{f}\|_{E_q}^{\frac{1}{p}}
\end{equation}
for constants $c=c(d,p,l)$ and $\gamma=\gamma(q,p)>0$. (Actually,  $\mathsf{f}$ in Corollary \ref{cor3} satisfies \eqref{H}, i.e.\,$\mathsf{f}=\mathsf{f}_{\mathfrak s}+\mathsf{f}_{\mathfrak b}$ where $|\mathsf{f}_{\mathfrak s}| \in E_q$, $q>1$ and $\mathsf{f}_{\mathfrak b}$ is bounded on $\mathbb R^{d+1}$. The term $\mathsf{f}_{\mathfrak b}$ is dealt with using the fact that $\rho \in L^p(\mathbb R^{d})$, so we only need to apply \eqref{rho_fin} to $|\mathsf{f}_{\mathfrak s}|$.)

To prove \eqref{rho_fin}, we estimate
\begin{align*}
& \|\rho_y\mathbf{1}_{[r,T]} |\mathsf{f}|^{\frac{1}{p}}\|_p^p  =\langle \rho_y^p \mathbf{1}_{[r,T]}|\mathsf{f}|\rangle \leq \sum_{k=0}^\infty (1+lk^2)^{-\nu p} \langle \mathbf{1}_{[r,T]}|\mathsf{f}| \mathbf{1}_{C_{k+1}(r,y)-C_k(r,y)}\rangle \qquad C_0(r,y):=\varnothing\\
& \leq \sum_{k=1}^\infty  (1+l k^2)^{-\nu p} \langle  \mathbf{1}_{[r,T]}|\mathsf{f}| \mathbf{1}_{C_{k+1}(r,y)}\rangle \\
& \leq \sum_{k=1}^\infty  (1+l k^2)^{-\nu p} |C_{k+1}(r,y) \cap ([r,T] \times \mathbb R^d) |^{\frac{1}{q'}}\frac{|C_{k+1}|^{\frac{1}{q}}}{k+1}(k+1)\biggl(\frac{1}{|C_{k+1}|}\langle |\mathsf{f}|^q \mathbf{1}_{C_{k+1}(r,y)}\rangle \biggr)^{\frac{1}{q}} \\
& \leq c_0(d)\sum_{k=1}^\infty  (1+l k^2)^{-\nu p} |T-r|^{\frac{1}{q'}}(k+1)^\frac{d}{q'}\frac{|C_{k+1}|^{\frac{1}{q}}}{k+1}\|\mathsf{f}\|_{E_q}\\
& \leq c_1(d,l)|T-r|^{\frac{1}{q'}}\sum_{k=1}^\infty k^{-2\nu p}k^{\frac{d}{q'}}k^{\frac{d+2}{q}-1}\|\mathsf{f}\|_{E_q} =:c^p|T-r|^{\frac{1}{q'}}\|\mathsf{f}\|_{E_q},
\end{align*}
where $c=c(d,p,l)<\infty$ since, by our choice of $\nu$ in \eqref{eta_two_est}, $-2\nu p + \frac{d}{q'}+\frac{d+2}{q}-1<-1$.   \hfill \qed

\bigskip

\section{Proof of Theorem \ref{thm3}}
\label{proof_thm3_sect}

Below we follow an argument from \cite{KiM} but use different embeddings, i.e.\,the ones established in Corollaries \ref{cor1}-\ref{cor3}.

The first assertion (\textit{i}), i.e.\,for all $x \in \mathbb R^d$, $0 \leq t \leq r \leq T$,
$$
\langle P^{t,r}(x,\cdot)\rangle=1,
$$
follows from \eqref{weight_bd0} with $\mathsf{f}=0$ and Theorem \ref{thm2}(\textit{i}). Namely, we first show for a fixed $x$, using weight $\rho$, that for every $\varepsilon>0$ there exists $R>0$ such that $\langle P_m^{t,r}(x,\cdot)\mathbf{1}_{\mathbb R^d - B(0,R)}(\cdot)\rangle<\varepsilon$ for all $m=1,2,\dots$, and so $\langle P_m^{t,r}(x,\cdot)\mathbf{1}_{B(0,R)}(\cdot)\rangle \geq 1- \varepsilon$. Passing to the limit in $m$ and then in $R \rightarrow \infty$, we obtain $\langle P^{t,r}(x,\cdot)\rangle \geq 1- \varepsilon$, which yields the required.

\medskip

(\textit{ii}) 
For every $n=1,2,\dots$, let $X^n_t=X^n_{t,x}$ denote the strong solution to the approximating SDE
$$
X^n_t=x-\int_0^t b_n(s,X^n_s)ds + \sqrt{2}B_t, \quad x \in \mathbb R^d,
$$
on a complete probability space ($\Omega$, $\mathcal F_t$, $\mathbf{P}$), where $\{b_n\}$ are given by \eqref{b_n}.

Step 1: There exists a constant $C>0$ independent of $n$, $k$ such that
\begin{equation}
\label{est0}
\sup_n \sup_{x \in \mathbb R^d}\mathbf E\int_{s}^r |b_k(t,X^n_{t,x})|dt \leq CF(r-s)
\end{equation}
for $0 \leq s \leq r \leq T$, where $$F(h):=h^\gamma,$$ constants $C$ and $\gamma>0$ (from Corollary \ref{cor3}) are independent of $n$ and $k$.
Indeed, let $v=v_{n,k}$ be the solution to the terminal-value problem
$$
(\partial_t + \Delta - b_n \cdot \nabla)v=-|b_k|, \quad v(r,\cdot)=0, \quad t \leq r.
$$
By It\^{o}'s formula,
$$
v(r,X_r^n)=v(s,X_s^n) + \int_s^r (\partial_t v + \Delta v - b_n \cdot \nabla v)(t,X^n_t) dt + \sqrt{2}\int_s^r \nabla v(t,X_t^n)dB_t,
$$
hence
$$
0=v(s,X_s^n) - \int_s^r |b_k(t,X_t^n)| dt + \sqrt{2}\int_s^r \nabla v(t,X_t^n)dB_t.
$$
Taking expectation, we obtain
$$
\mathbf E\int_s^r |b_n(t,X_t^n)|dt = \mathbf E v(s,X_s^n).
$$
Since $\mathbf E v(s,X_s^n) \leq \|v(s,\cdot)\|_{L^\infty(\mathbb R^d)}$, we obtain 
from Corollary \ref{cor3} with $\mathsf{f}_k=\mathsf{b}_k$ and initial data $g=0$
\begin{align*}
\mathbf E\int_s^r |b_k(t,X_t^n)|dt  \leq C(r-s)^\gamma
\end{align*}
with constants $C$ and $\gamma>0$ independent of $k$, $n$ $\Rightarrow$ \eqref{est0}.

\medskip

By a standard result (see e.g.\,\cite[Ch.\,2]{GvC}), given a conservative backward Feller evolution family, there exist probability measures $\mathbb P_x$ $(x \in \mathbb R^d)$ on
$(D([0,T],\mathbb R^d),\mathcal B'_t=\sigma(\omega_r \mid 0\leq r \leq t))$, where $D([0,T],\mathbb R^d)$ is the space of right-continuous functions having left limits, and $\omega_t$ is the coordinate process, such that
$$
\mathbb E_{x}[f(\omega_r)]=P^{0,r}f(x), \quad 0 \leq r \leq T. 
$$
Here and below, $\mathbb E_x:=\mathbb E_{\mathbb P_x}$. Also, put $\{\mathbb P_x^n:=(\mathbf PX^n)^{-1}\}_{n=1}^\infty$ and set $\mathbb E_x^n:=\mathbb E_{\mathbb P_x^n}$.

\medskip

Step 2: $\mathbb E_x[\int_0^r |b(t,\omega_t)|dt]<\infty$.

\smallskip

Indeed, by Step 1, $\sup_n \sup_{x \in \mathbb R^d}\mathbb E_x^n\int_{s}^r |b_k(t,\omega_t)|dt \leq CF(r-s)$.
Hence by the convergence result in Corollary \ref{cor2} (with $\mathsf{f}:=\mathbf{1}_{B(0,k)}b_k$) 
$$\mathbb E_x[\int_s^r |\mathbf{1}_{B(0,k)}b_k(t,\omega_t)|dt] \leq CF(r-s)<\infty.
$$ It remains to apply Fatou's Lemma in $k$.

\medskip

Step 3: For every $f \in C_c^2(\mathbb R^d)$, the process
\begin{equation}
\label{M}
M_r^f:= f(\omega_r)-f(x) + \int_0^r (-\Delta f + b \cdot \nabla f)(t,\omega_t)dt
\end{equation}
is a $\mathcal B'_r$-martingale under $\mathbb P_x$. 

Indeed, let us note first that
\begin{equation}
\label{conv0}
\tag{$\star$}
\mathbb E^m_x[f(\omega_r)] \rightarrow \mathbb E_x[f(\omega_r)], \quad \mathbb E^m_x[\int_0^r (-\Delta f)(\omega_t)dt] \rightarrow \mathbb E_x[\int_0^r (-\Delta f)(\omega_t)dt] \quad (m \rightarrow \infty),
\end{equation}
as follows from the convergence result in Theorem \ref{thm2}(\textit{i}).
Next, we note that
\begin{equation}
\label{conv1}
\tag{$\star\star$}
\mathbb E^m_x \int_0^r (b_m \cdot \nabla f)(t,\omega_t)dt \rightarrow \mathbb E_x \int_0^r (b \cdot \nabla f)(t,\omega_t)dt \quad (m \rightarrow \infty).
\end{equation}
The latter follows from
\begin{align}
\mathbb E^m_x\bigg|\int_0^r \big((b_m-b_n)\cdot \nabla f\big)(t,\omega_t)dt\, \eta(\omega) \bigg| & \leq C\|\eta\|_\infty \|\mathbf{1}_{[0,r]}|b_m-b_n|^{\frac{1}{p}}|\nabla f|\|_p
\label{a}
\tag{$a$}
\end{align}
as $m,n \rightarrow \infty$;
\begin{equation}
\label{b}
\tag{$b$}
\mathbb E^m_x\bigg[\int_0^r (b_n \cdot \nabla f)(t,\omega_t)dt  \cdot \eta(\omega) \bigg] \rightarrow \mathbb E_x\bigg[\int_0^r (b_n \cdot \nabla f)(t,\omega_t)dt \cdot \eta(\omega) \bigg]
\end{equation}
as $m \rightarrow \infty$;
\begin{align}
\label{c}
\tag{$c$}
\mathbb E_x\bigg|\int_0^r \big((b-b_n)\cdot \nabla f\big)(t,\omega_t)dt\, \eta(\omega) \bigg| & \leq C\|\eta\|_\infty\|\mathbf{1}_{[0,r]}|b-b_n|^{\frac{1}{p}}|\nabla f|\|_p \rightarrow 0
\end{align}
as $n \rightarrow \infty$. 
The proof of the inequality in (a) follows the proof of \eqref{est0} but uses Corollary \ref{cor2} with $\mathsf{f}=b_n-b_m$ and $g=0$. The convergence in (a) follows from the fact that $b_n-b_m \rightarrow 0$ in $[L^1_{\loc}(\mathbb R^{d+1})]^d$.
Assertion (b) follows from Corollary \ref{cor2}. 
The proof of (c) is similar to the proof of (a) except that we pass to the limit in $m$ and then in $k$ using Fatou's Lemma.

Now, since
$$
M_{r,m}^{f}:=f(\omega_r)-f(x) + \int_0^r (-\Delta f + b_m \cdot \nabla f)(t,\omega_t)dt
$$
is a $\mathcal B'_r$-martingale under $\mathbb P^m_x$, 
$$
x \mapsto \mathbb E^m_x[f(\omega_r)] - f(x) +\mathbb E^m_x\int_0^r (-\Delta f + b_m\cdot\nabla f)(t,\omega_t)dt \quad \text{ is identically zero on } \mathbb R^d,
$$
and so by \eqref{conv0}, \eqref{conv1}
$$
x \mapsto \mathbb E_x[f(\omega_r)] - f(x) +\mathbb E_x\int_0^r (-\Delta f + b\cdot\nabla f)(t,\omega_t)dt \quad  \text{ is identically zero in } \mathbb R^d.
$$
Since $\{\mathbb P_x\}_{x \in \mathbb R^d}$ are determined by a Feller evolution family, and thus constitute a Markov process, we can conclude (see e.g.\,the proof of \cite[Lemma 2.2]{Kr1}) that $M_r^f$ is a $\mathcal B'_r$-martingale under $\mathbb P_x$.

\medskip

Step 4:
$\{\mathbb P_x\}_{x \in \mathbb R^d}$ are concentrated on 
$(C([0,T],\mathbb R^d),\mathcal B_t)$.

\smallskip

By Step 3, $\omega_t$ is a semimartingale under $\mathbb P_x$, so It\^{o}'s formula yields, for every $g \in C_c^\infty(\mathbb R^d)$, that
\begin{equation}
\label{g_}
g(\omega_t)-g(x)=\sum_{s \leq t}\bigl(g(\omega_s)-g(\omega_{s-})\bigr) + S_t,
\end{equation}
where $S_t$ is defined in terms of some  integrals and sums of $(\partial_{x_i}g)(\omega_{s-})$ and $(\partial_{x_i}\partial_{x_j}g)(\omega_{s-})$ in $s$, see \cite[Sect.\,2]{CKS} for details.
Now,
let $A$, $B$ be arbitrary compact sets in $\mathbb R^d$ such that $\dist(A,B)>0$. 
Fix $g \in C_c^\infty(\mathbb R^d)$ that separates $A$, $B$, say, $g = 0$ on $A$, $g = 1$ on $B$. Set
$$
K^g_t:=\int_0^t \mathbf{1}_A(\omega_{s-})dM_s.
$$
In view of \eqref{M} and \eqref{g_}, when evaluating $K^g_t$ one needs to integrate $\mathbf{1}_A(\omega_{s-})$ with respect to $S_t$, however, one obtains zero since $(\partial_{x_i}g)(\omega_{s-})=(\partial_{x_i}\partial_{x_j}g)(\omega_{s-})=0$ if $\omega_{s-} \in A$. Thus,
\begin{align*}
K^g_t &=\sum_{s \leq t} \mathbf{1}_A \left(\omega_{s-}\right)g(\omega_{s}) +
\int_0^t \mathbf{1}_A(\omega_{s-})\bigl(-\Delta g + b\cdot\nabla g \bigr)(\omega_s)ds \\
&=\sum_{s \leq t} \mathbf{1}_A \left(\omega_{s-}\right)g(\omega_s).
\end{align*}
Since $M^g_t$ is a martingale, so is $K^g_t$. Thus, $\mathbb{E}_x\bigl[\sum_{s \leq t} \mathbf{1}_A (\omega_{s-})g(\omega_s)\bigr]=0.$ Using the Dominated Convergence Theorem, we further obtain
$\mathbb{E}_x\bigl[\sum_{s \leq t} \mathbf{1}_A (\omega_{s-})\mathbf{1}_B(\omega_s)\bigr]=0$, which yields the required. (By the way, this construction, in a more general form, can be used to control the jumps of stable process perturbed by a drift, see \cite{CKS}.)

\medskip

We denote the restriction of $\mathbb P_x$ from $(D([0,T],\mathbb R^d), \mathcal B_t')$  to $(C([0,T],\mathbb R^d),\mathcal B_t)$ again by $\mathbb P_x$, and thus obtain that
for every $x \in \mathbb R^d$ and all $f \in C_c^2(\mathbb R^d)$ 
$$
M_r^f=f(\omega_r) - f(x) + \int_0^r (-\Delta f + b\cdot\nabla f)(t,\omega_t)dt, \quad \omega \in C([0,T],\mathbb R^d),$$
is a $\mathcal B_r$-martingale under $\mathbb P_x$.

Thus, $\mathbb P_x$ is a $\mathcal B_r$-martingale solution to \eqref{seq}. (Alternatively, we could have used a tightness argument, cf.\,\cite{KiM}.)

To show that $\mathbb P_x$ is a weak solution it suffices to show that $M_r^f$ is also a martingale for $f(x)=x_i$ and $f(x)=x_ix_j$, which is done by following closely \cite[proof of Lemma 6]{KiS1} and employing weight $\rho$ and \eqref{weight_bd0} in Corollary \ref{cor3}.

\medskip

(\textit{iii}) This follows from Corollary \ref{cor2}, cf.\,proof of (\textit{ii}) above. 

($\textit{iii}'$) is proved by retracing the steps in the proof of (\textit{iii}) but now using the Sobolev embedding property of the operator $(\lambda+\partial_t-\Delta)^{-\frac{1}{2p'}}$ in the solution representation of Theorem \ref{thm1}, i.e.\,the fact that $(\lambda+\partial_t-\Delta)^{-\frac{1}{2p'}}$ bounded as an operator from $L^{\nu}$, $\nu=\frac{p(d+2)}{d+p+1}$ ($<p$) to $L^p$. Note that by selecting $p>d+1$ close to $d+1$ we make $\nu>\frac{d+2}{2}$ close to $\frac{d+2}{2}$.

\medskip

(\textit{iv}) Suppose that there exist $\mathbb P^1_x$, $\mathbb P^2_x$, two martingale solutions to \eqref{seq}, that satisfy
\begin{align}
\label{kr_est2}
\mathbb E_{x}^i\int_0^T |b(r,\omega_t)h(t,\omega_t)|dt  \leq c\|\mathbf{1}_{[0,T]}b^{\frac{1}{p}}h\|_p, \quad  h \in C_c(\mathbb R^{d+1})
\end{align}
with constant $c$ independent of $h$  ($i=1,2$).
Here and below, $\mathbb E_{x}^1:=\mathbb E_{\mathbb P_x^1}$, $\mathbb E_{x}^2:=\mathbb E_{\mathbb P_x^2}$.
We will show that for every $F \in C_c(\mathbb R^{d+1})$ we have
\begin{equation}
\label{id4}
\mathbb E_x^1[\int_0^T F(t,\omega_t)dt]=\mathbb E_x^2[\int_0^T F(t,\omega_t)dt],
\end{equation}
which implies $\mathbb P_x^1=\mathbb P_x^2$.

Proof of \eqref{id4}. Let $u_n \in C([0,T],C_\infty(\mathbb R^d))$ be the classical solution to 
\begin{equation}
\label{eq_F}
(\partial_t + \Delta + b_n \cdot \nabla)u_n=F, \quad u_n(T,\cdot)=0,
\end{equation}
where, recall, $b_n=\mathbf{1}_n b$, and $\mathbf{1}_n$ is the indicator of $\{|b| \leq n\}$.
Set $\tau_R:=\inf\{t \geq 0 \mid |\omega_t| \geq R\}$, $R>0$.
By It\^{o}'s formula 
\begin{align}
\mathbb E_x^i u_n(T \wedge \tau_R,\omega_{T \wedge \tau_R})  & = u_n(0,x)+\mathbb E_x^i \int_0^{T \wedge \tau_R} F(t,\omega_t)dt \notag \\
& +  \mathbb E_x^i \int_0^{T \wedge \tau_R} \big[(b-b_n)\cdot \nabla u_n\big](t,\omega_t)dt \label{e_i}
\end{align}
($i=1,2$).
We have
\begin{align}
\biggl| \mathbb E_x^i \int_0^{T \wedge \tau_R} \big[(b-b_n)\cdot \nabla u_n\big](t,\omega_t)dt \biggr| & \leq \mathbb E_x^i \int_0^{T \wedge \tau_R} \big[|b|(1-\mathbf{1}_n) |\nabla u_n|\big](t,\omega_t)dt \label{bd_est} \\
& (\text{we are applying \eqref{kr_est2}}) \notag \\
& \leq c\|\mathbf{1}_{[0,T] \times B(0,R)}|b|^{\frac{1}{p}}(1-\mathbf{1}_n)|\nabla u_n|\|_p. \notag
\end{align}
At this point we note that $\tilde{u}_n(t):=e^{\lambda (T-t)}u_n(t)$ satisfies 
$$
(\lambda + \partial_t + \Delta + b_n \cdot \nabla)u_n=\mathbf{1}_{[0,T]}e^{\lambda (T-t)}F.
$$
Hence we can apply to $|b|^{\frac{1}{p}}|\nabla u_n|$ the solution representation of Corollary \ref{cor2} (after reversing time and taking there $g=0$). Using $|b| \geq |b_n|$, we then obtain an independent on $n$ $L^p(\mathbb R^{d+1})$ majorant on $|b|^{\frac{1}{p}}|\nabla u_n|$. Therefore, since $1-\mathbf{1}_n \rightarrow 0$ a.e.\,on $\mathbb R^{d+1}$ as $n \rightarrow \infty$, we have  
$$
\mathbb E_x^i \int_0^{T \wedge \tau_R} \big[(b-b_n)\cdot \nabla u_n\big](t,\omega_t)dt \rightarrow 0 \quad (n \rightarrow \infty).
$$

We are left to note, using again Corollary \ref{cor2}, that solutions $u_n$ converge to a function $u \in C([0,T], C_\infty(\mathbb R^d))$. Therefore, we can pass to the limit in \eqref{e_i}, first in $n$ and then in $R \rightarrow \infty$, to obtain
$$
0 =u(0,x)+\mathbb E_x^i \int_0^{T} F(t,\omega_t)dt \quad i=1,2,
$$
which gives \eqref{id4}.

\medskip

($\textit{iv}'$) is proved by following the proof of (\textit{iv}) but estimating \eqref{bd_est} differently. 
By our assumption, we have
\begin{align*}
\mathbb E_{x}^i\int_0^T |h(t,\omega_t)|dt  \leq c\|\mathbf{1}_{[0,T]}h\|_\nu, \quad  h \in C_c(\mathbb R^{d+1}), \quad i=1,2,
\end{align*}
where, recall, $\nu>\frac{d+2}{2}$ is close to $\frac{d+2}{2}$.
Therefore,
\begin{align*}
\mathbb E_x^i \int_0^{T \wedge \tau_R} \big[|b|(1-\mathbf{1}_n) |\nabla u_n|\big](t,\omega_t)dt & \leq c\|\mathbf{1}_{[0,T] \times B(0,R)}|b|(1-\mathbf{1}_n)|\nabla u_n|\|_\nu \\
& \leq
 c\|\mathbf{1}_{[0,T] \times B(0,R)}|b|(1-\mathbf{1}_n)\|_{s'}\|\nabla u_n\|_s, \qquad \frac{1}{s}+\frac{1}{s'}=\frac{1}{\nu},
\end{align*}
It follows from the solution representation of Corollary \ref{cor2} and the parabolic Sobolev embedding theorem that
$$
\|\nabla u_n\|_s \leq C\|\mathbf{1}_{[0,T]}e^{\lambda T}F\|_{p} \quad \text{ for } s=\frac{d+2}{d+1}p\;\; \text{ close to } \frac{d+2}{d+1}p.
$$
Assuming that the Morrey norm $\|b_{\mathfrak s}\|_{E_q}$ is sufficiently small, we can select $p$ sufficiently large to make $s'>\nu$ close to $\nu$ and hence close to $\frac{d+2}{2}$. (That is, since by our assumption $|b| \in L_{\loc}^{\frac{d+2}{2}+\varepsilon}$ for some $\varepsilon>0$, we need $s' \geq \frac{d+2}{2}+\varepsilon$.) Now, since $1-\mathbf{1}_n \rightarrow 0$ a.e.\,on $\mathbb R^{d+1}$ as $n \rightarrow \infty$, we have  $\|\mathbf{1}_{[0,T] \times B(0,R)}|b|(1-\mathbf{1}_n)\|_{s'} \rightarrow 0$ as $n \rightarrow \infty$. The rest repeats the proof of (\textit{iv}). \hfill \qed

\end{document}